\documentclass[11pt,reqno]{amsart}
\usepackage{amsfonts}
\usepackage{amssymb}
\usepackage{amsmath}
\usepackage{eepic}
\usepackage{dsfont}

\allowdisplaybreaks

\newtheorem{theorem}{Theorem}[section]
\newtheorem{lemma}[theorem]{Lemma}
\newtheorem{proposition}[theorem]{Proposition}
\newtheorem{corollary}[theorem]{Corollary}

\theoremstyle{definition}

\newtheorem{definition}[theorem]{Definition}

\theoremstyle{plain}

\numberwithin{equation}{section}

\newcommand{\ga}{\gamma}
\newcommand{\Ac}{\mathcal A}

\newcommand{\Dc}{\mathcal D}
\newcommand{\GG}{\mathds G}
\newcommand{\Ec}{\mathcal E}
\newcommand{\Fc}{\mathcal F}
\newcommand{\Sc}{\mathds S}
\newcommand{\Gc}{\mathcal G}

\newcommand{\Ga}{\Gamma}
\newcommand{\inv}{^{-1}}

\def\arr#1{\mathop{{\buildrel #1\over\longrightarrow}}}
\def\ov{\overline}
\def\wh{\widehat}
\def\St#1{\mathrm{Stab}(#1)}

\begin{document}
    \title[Automata over a binary alphabet generating free groups]{Automata over a binary
      alphabet generating free groups of
      even rank}

\author{Benjamin Steinberg}\address{School of Mathematics and
  Statistics\\ Carleton University\\ 1125 Colonel By Drive\\ Ottawa\\
  Ontario  K1S 5B6 \\ Canada}
\email{bsteinbg@math.carleton.ca}
\thanks{The first author acknowledges the support of NSERC}
\author{Mariya Vorobets}\address{Department of Mathematics \\
Texas A\&M University\\ College Station, TX 77843-3368}
\email{mvorobet@math.tamu.edu}
\author{Yaroslav Vorobets}\address{Department of Mathematics \\
Texas A\&M University\\ College Station, TX 77843-3368}
\email{yvorobet@math.tamu.edu}
\thanks{The third author is supported by a Clay Research Scholarship}
\date{Version of \today}

\begin{abstract}
We construct automata over a binary alphabet with $2n$ states,
$n\geq 2$, whose states freely generate a free group of rank $2n$.
Combined with previous work, this shows that a free group of every
finite rank can be generated by finite automata over a binary
alphabet.  We also construct free products of cyclic groups of order
two via such automata.
\end{abstract}
\maketitle 

\section{Introduction}
The binary odometer is a finite state automaton over a binary
alphabet generating a free group of rank $1$. For a long time, it
was an open question whether there is a finite state automaton
generating a non-abelian free group~\cite{GNS}.  The first examples
were provided by Glasner and Mozes~\cite{glasner}.  The smallest
examples they obtained were a $14$-state automaton over a $6$-letter
alphabet generating a free group of rank $7$ and a $6$-state
automaton over a $14$-letter alphabet generating a free group of
rank $3$.
  Nekrashevych afterwards
constructed an automaton with $6$ states over a binary alphabet
generating a free group of rank two~\cite{selfsimilar}; again the
states of the automaton are not free generators and, in fact,
it was shown in \cite{GNS} that no two state automaton over a
binary alphabet generates a free group.

   Brunner and
Sidki conjectured in~\cite{Brunner} that the states of a certain
$3$-state automaton over a binary alphabet, constructed by
Aleshin~\cite{aleshin} in 1983, freely generate a free group of rank
$3$. Aleshin constructed this automaton, together with a certain
five state automaton, in order to construct a free subgroup of rank
$2$ in the group of all finite state transformations, but his
proof~\cite{aleshin} is not complete.  In~\cite{free1} the second
and third authors answered positively the conjecture of Brunner and
Sidki, showing that the Aleshin automaton freely generates a free
group of rank $3$.  In a sequel paper~\cite{free2}, they considered
a series of Aleshin-type automata with $2n+1$ states, $n\geq 1$,
over a binary alphabet that freely generate a free group of rank
$2n+1$.  Moreover, it was shown that the disjoint union of any two
distinct automata in this series generates their free product.  Thus
they established that, for every possible rank $n\geq 3$, except $4$
and $6$, there is a free group of rank $n$ generated by an $n$-state
automaton over a binary alphabet.  In the process, Aleshin's
claim~\cite{aleshin} was verified.  It should be noted, however, that
connected automata were constructed only for free groups of odd rank.

In this paper we construct a new family of finite automata over a
binary alphabet generating free groups.  Our family contains, for
any $n\geq 2$, a connected $2n$-state automaton freely generating a
free group of rank $2n$.  In particular, combined with the results
of~\cite{free1,free2} and~\cite{selfsimilar}, this shows that a free
group of any rank can be generated by a finite state automaton over
a binary alphabet and, moreover, if $n\geq 3$, then an $n$-state
connected automaton does the job.

Nekrashevych's construction~\cite{selfsimilar} was based on a result
of Muntyan and Savchuk~\cite[Theorem 1.4]{selfsimilar}, showing that
a certain $3$-state automaton over a binary alphabet, related to
Aleshin's automaton, generates a free product of three cyclic groups
of order two.   The second and third authors, in~\cite{free2}, in
turn showed that there are connected $2n+1$ state automata over a
binary alphabet generating a free product of $2n+1$ cyclic groups of
order two for any $n\geq 1$.  Again the disjoint union of distinct
automata from this family generate their free product and so they
were able to construct, for any $n\geq 3$, a finite state automaton
generating a free product of $n$ cyclic groups of order two (in this
setting they could add a single isolated state to obtain products of
$4$ and $6$ cyclic groups of order two). In this paper we use our
family to construct, for any $n\geq 2$, a $2n$-state connected
automaton over a binary alphabet generating a free product of $2n$
cyclic groups of order two.

All automata discussed in this paper are bireversible and it is an
open problem to construct an automaton that is not bireversible
generating a free non-abelian group. By an unpublished result of
Abert, such an automaton cannot be contracting in the sense
of~\cite{selfsimilar}.

The paper is organized as follows.  The second section gives the
reader the basic background about groups generated by finite state
automata, with a special emphasis on dual automata and bireversible
automata. In particular, we discuss our conventions for dealing with
dual automata.  The third section introduces the automata that will
play a key role in this paper.  The fourth section proves the main
results on freeness.

\section{Automaton groups}\label{sec:automatongroups}
In this section, we collect some of the basic notions from the
theory of groups generated by finite state automata, also called
automaton groups.  This is a special case of the notion of a
self-similar group~\cite{selfsimilar}, due to Nekrashevych, but is
the principal case that has been studied.   For more information
consult~\cite{GNS,selfsimilar}.

\subsection{Preliminaries and notation for free monoids}
  First some preliminaries and notation.  If $A$ is a
finite alphabet (that is a finite set), then $A^*$ denotes the free
monoid on $A$.  If $w=a_1\cdots a_n\in A^*$, then the
\emph{reversal}  of $w$ is the word $w^{\rho} = a_n\cdots a_1$.  We
set $A^{\pm} = A\cup A\inv$, where $A\inv$ is a disjoint set in
bijection with $a$ via a map $a\mapsto a\inv$. Then $(A^{\pm})^*$ is
the \emph{free monoid with involution}, where the involution is
defined in the usual way. The group of involution-preserving
automorphisms of $(A^{\pm})^*$ is isomorphic to $S_A$ via the map
sending $\sigma\in S_A$ to the map $\sigma^{\pm}:(A^{\pm})^*\to
(A^{\pm})^*$ defined by
\[\sigma^{\pm}(a_1^{e_1}\cdots a_n^{e_n}) =
\sigma(a_1)^{e_1}\cdots \sigma (a_n)^{e_n}\] with $a_i\in A$ and
$e_i\in \{\pm 1\}$.  Sometimes we shall wish to distinguish
$\sigma^{\pm}$ from the permutation of $A^{\pm}$ obtained by
restricting $\sigma^{\pm}$ to letters.  So if $\sigma\in S_A$, then
$\ov \sigma\in S_{A^{\pm}}$ denotes the permutation defined by $\ov
\sigma(a^e)=\sigma(a)^{e}$ for $a\in A$, $e=\pm 1$.  The map
$\sigma\mapsto \ov \sigma$ is of course a monomorphism $S_A\to
A_{A^{\pm}}$. To each element $\sigma\in S_A$, we can also associate
an automorphism $\sigma^*:A^*\to A^*$ by defining
$\sigma^*(a)=\sigma(a)$.  Note that $\sigma^\pm = \ov{\sigma}^*$.

\subsection{Mealy automata}
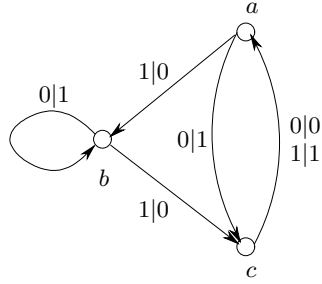
\begin{figure}[htbp]
\begin{center}
\setlength{\unitlength}{0.00062500in}
\begingroup\makeatletter\ifx\SetFigFont\undefined%
\gdef\SetFigFont#1#2#3#4#5{%
  \reset@font\fontsize{#1}{#2pt}%
  \fontfamily{#3}\fontseries{#4}\fontshape{#5}%
  \selectfont}%
\fi\endgroup%
\begin{picture}(2582,2395)(0,-10)
\put(453.000,1171.000){\arc{3632.267}{5.7634}{6.7690}}
\blacken\path(2111.592,1980.921)(2029.000,2073.000)(2058.544,1952.887)(2068.247,1998.733)(2111.592,1980.921)
\put(3515.031,1245.275){\arc{3632.065}{2.6216}{3.6273}}
\blacken\path(1856.406,435.077)(1939.000,343.000)(1909.453,463.112)(1899.751,417.266)(1856.406,435.077)
\put(784,1200){\ellipse{150}{150}}
\put(1984,2100){\ellipse{150}{150}}
\put(1984,300){\ellipse{150}{150}}
\path(856,1152)(1921,327)
\blacken\path(1807.762,376.771)(1921.000,327.000)(1844.506,424.204)(1854.594,378.441)(1807.762,376.771)
\path(1893,2075)(853,1227)
\blacken\path(927.044,1326.083)(853.000,1227.000)(964.960,1279.582)(918.102,1280.083)(927.044,1326.083)
\path(726,1246)(725,1247)(721,1251)
    (716,1256)(707,1264)(695,1275)
    (681,1288)(664,1303)(645,1318)
    (625,1335)(603,1351)(580,1367)
    (556,1382)(531,1397)(504,1409)
    (475,1420)(445,1429)(412,1436)
    (378,1439)(342,1439)(310,1435)
    (280,1429)(252,1421)(227,1413)
    (205,1404)(185,1396)(168,1387)
    (153,1379)(140,1372)(128,1364)
    (116,1356)(105,1349)(94,1340)
    (84,1331)(72,1321)(61,1309)
    (50,1295)(39,1280)(28,1262)
    (20,1243)(14,1221)(12,1199)
    (15,1177)(22,1155)(32,1135)
    (44,1117)(56,1100)(68,1086)
    (80,1073)(92,1062)(104,1052)
    (115,1042)(127,1033)(139,1025)
    (152,1016)(166,1008)(181,999)
    (199,989)(219,980)(241,970)
    (267,961)(295,953)(325,946)
    (357,942)(392,942)(426,945)
    (457,952)(487,962)(514,974)
    (539,987)(563,1002)(585,1018)
    (606,1036)(626,1053)(645,1071)
    (662,1088)(677,1103)(690,1117)
    (701,1129)(718,1148)
\blacken\path(660.342,1038.567)(718.000,1148.000)(615.627,1078.575)(661.989,1085.400)(660.342,1038.567)
\put(1984,50){\makebox(0,0)[lb]{{\SetFigFont{9}{10.8}{\rmdefault}{\mddefault}{\updefault}$c$}}}
\put(754,810){\makebox(0,0)[lb]{{\SetFigFont{9}{10.8}{\rmdefault}{\mddefault}{\updefault}$b$}}}
\put(1984,2260){\makebox(0,0)[lb]{{\SetFigFont{9}{10.8}{\rmdefault}{\mddefault}{\updefault}$a$}}}
\put(252,1467){\makebox(0,0)[lb]{{\SetFigFont{9}{10.8}{\rmdefault}{\mddefault}{\updefault}$0|1$}}}
\put(1082,1655){\makebox(0,0)[lb]{{\SetFigFont{9}{10.8}{\rmdefault}{\mddefault}{\updefault}$1|0$}}}
\put(1082,505){\makebox(0,0)[lb]{{\SetFigFont{9}{10.8}{\rmdefault}{\mddefault}{\updefault}$1|0$}}}
\put(1412,1085){\makebox(0,0)[lb]{{\SetFigFont{9}{10.8}{\rmdefault}{\mddefault}{\updefault}$0|1$}}}
\put(2352,1205){\makebox(0,0)[lb]{{\SetFigFont{9}{10.8}{\rmdefault}{\mddefault}{\updefault}$0|0$}}}
\put(2352,965){\makebox(0,0)[lb]{{\SetFigFont{9}{10.8}{\rmdefault}{\mddefault}{\updefault}$1|1$}}}
\end{picture}

\end{center}
\caption{Aleshin's automaton $\mathsf A$\label{aleshindiagram}}
\end{figure}

A \emph{finite (Mealy) automaton} \cite{Eilenberg,Arbib} $\Ac$ is a
$4$-tuple $(Q,A,\delta,\lambda)$ where $Q$ is a finite set of
states, $A$ is a finite alphabet, $\delta:Q\times A\to Q$ is the
transition function and $\lambda: Q\times A\to A$ is the output
function.  We shall always write, for $q\in Q$ and $a\in A$, $\delta
(q,a) = q_a$ and $\lambda (q,a) = q(a)$.  Automata are usually
represented by so-called \emph{Moore diagrams}.  The Moore diagram
for $\mathcal A$ is a directed graph with vertex set $Q$. The edges
are of the form $q\arr{a\mid q(a)} q_a$.  For example Aleshin's
automaton~\cite{aleshin} is given by the Moore diagram in
Figure~\ref{aleshindiagram}.

For instance, in state $b$ with input $1$, the automaton outputs $0$
and moves to state $c$. Sometimes we shall just draw the transitions
and omit the output from the Moore diagram;  the resulting graph is
called the \emph{transition diagram}.  The transition and output
functions of an automaton $\Ac$ extend inductively to the free
monoid $A^*$ via the rules:
\begin{gather} q_{au} = (q_a)_u
  \\\label{action} q(au) = q(a)q_a(u)
\end{gather}
where $a\in A$, $u\in A^*$.

\begin{figure}[htbp]
\begin{center}
\setlength{\unitlength}{0.00062500in}
\begingroup\makeatletter\ifx\SetFigFont\undefined%
\gdef\SetFigFont#1#2#3#4#5{%
  \reset@font\fontsize{#1}{#2pt}%
  \fontfamily{#3}\fontseries{#4}\fontshape{#5}%
  \selectfont}%
\fi\endgroup%
\begin{picture}(2582,2395)(0,-10)
\put(453.000,1171.000){\arc{3632.267}{5.7634}{6.7690}}
\blacken\path(2111.592,1980.921)(2029.000,2073.000)(2058.544,1952.887)(2068.247,1998.733)(2111.592,1980.921)
\put(3515.031,1245.275){\arc{3632.065}{2.6216}{3.6273}}
\blacken\path(1856.406,435.077)(1939.000,343.000)(1909.453,463.112)(1899.751,417.266)(1856.406,435.077)
\put(784,1200){\ellipse{150}{150}}
\put(1984,2100){\ellipse{150}{150}}
\put(1984,300){\ellipse{150}{150}}
\path(856,1152)(1921,327)
\blacken\path(1807.762,376.771)(1921.000,327.000)(1844.506,424.204)(1854.594,378.441)(1807.762,376.771)
\path(1893,2075)(853,1227)
\blacken\path(927.044,1326.083)(853.000,1227.000)(964.960,1279.582)(918.102,1280.083)(927.044,1326.083)
\path(726,1246)(725,1247)(721,1251)
    (716,1256)(707,1264)(695,1275)
    (681,1288)(664,1303)(645,1318)
    (625,1335)(603,1351)(580,1367)
    (556,1382)(531,1397)(504,1409)
    (475,1420)(445,1429)(412,1436)
    (378,1439)(342,1439)(310,1435)
    (280,1429)(252,1421)(227,1413)
    (205,1404)(185,1396)(168,1387)
    (153,1379)(140,1372)(128,1364)
    (116,1356)(105,1349)(94,1340)
    (84,1331)(72,1321)(61,1309)
    (50,1295)(39,1280)(28,1262)
    (20,1243)(14,1221)(12,1199)
    (15,1177)(22,1155)(32,1135)
    (44,1117)(56,1100)(68,1086)
    (80,1073)(92,1062)(104,1052)
    (115,1042)(127,1033)(139,1025)
    (152,1016)(166,1008)(181,999)
    (199,989)(219,980)(241,970)
    (267,961)(295,953)(325,946)
    (357,942)(392,942)(426,945)
    (457,952)(487,962)(514,974)
    (539,987)(563,1002)(585,1018)
    (606,1036)(626,1053)(645,1071)
    (662,1088)(677,1103)(690,1117)
    (701,1129)(718,1148)
\blacken\path(660.342,1038.567)(718.000,1148.000)(615.627,1078.575)(661.989,1085.400)(660.342,1038.567)
\put(1984,50){\makebox(0,0)[lb]{{\SetFigFont{9}{10.8}{\rmdefault}{\mddefault}{\updefault}$c\inv$}}}
\put(754,810){\makebox(0,0)[lb]{{\SetFigFont{9}{10.8}{\rmdefault}{\mddefault}{\updefault}$b\inv$}}}
\put(1984,2260){\makebox(0,0)[lb]{{\SetFigFont{9}{10.8}{\rmdefault}{\mddefault}{\updefault}$a\inv$}}}
\put(252,1467){\makebox(0,0)[lb]{{\SetFigFont{9}{10.8}{\rmdefault}{\mddefault}{\updefault}$1|0$}}}
\put(1082,1655){\makebox(0,0)[lb]{{\SetFigFont{9}{10.8}{\rmdefault}{\mddefault}{\updefault}$0|1$}}}
\put(1082,505){\makebox(0,0)[lb]{{\SetFigFont{9}{10.8}{\rmdefault}{\mddefault}{\updefault}$0|1$}}}
\put(1412,1085){\makebox(0,0)[lb]{{\SetFigFont{9}{10.8}{\rmdefault}{\mddefault}{\updefault}$1|0$}}}
\put(2352,1205){\makebox(0,0)[lb]{{\SetFigFont{9}{10.8}{\rmdefault}{\mddefault}{\updefault}$0|0$}}}
\put(2352,965){\makebox(0,0)[lb]{{\SetFigFont{9}{10.8}{\rmdefault}{\mddefault}{\updefault}$1|1$}}}
\end{picture}
\end{center}
\caption{The inverse of Aleshin's automaton\label{aleshininverse}}
\end{figure}
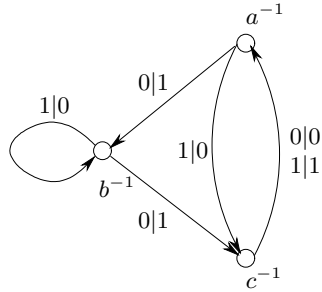

We use $\Ac_q$ to denote the \textit{initial automaton} $\Ac$ with
designated start state $q$. Sometimes, when no confusion can occur,
we use simply the letter $q$ to denote this initial automaton.
Abusing notation,  there is a function $\Ac_q:A^*\to A^*$ given by
$w\mapsto q(w)$.  For instance, with $\mathsf A$ the automaton
above, $\mathsf A_b(010) = 100$.    In general the function $\Ac_q$
is length preserving and preserves common prefixes.  It extends
continuously to the set of right infinite words $A^{\omega}$ via the
formula
\begin{equation}\label{actiononcantor}
\Ac_q(a_0a_1\cdots) = \lim_{n \to \infty} \Ac_q(a_0\cdots a_n)
\end{equation}
 where $A^{\omega}$
is given the product topology, making it homeomorphic to a Cantor
set~\cite{GNS}. If one defines a metric on $A^{\omega}$ by defining
$d(u,v) = |A|^{-|u\wedge v|}$, where $u\wedge v$ is the longest
common prefix of $u$ and $v$, then $\Ac_q$ is a metric contraction
(where we say a map $f$ is a \emph{metric contraction} if
$d(f(u),f(v))\leq d(u,v)$).  If, for each $q$, the function
$a\mapsto q(a)$ is a permutation, then each $\Ac_q$ is invertible
and induces an isometry of $A^{\omega}$~\cite{GNS,selfsimilar}. In
this case the automaton is called \emph{invertible}.  The
\emph{inverse automaton} $\Ac\inv$ of $\Ac$ is the automaton
obtained by taking the Moore diagram for $\Ac$ and swapping the left
hand side and right hand side of the edge labels.   We usually
denote the state of $\Ac\inv$ corresponding to $q$ by $q\inv$.   One
can check that $\Ac\inv_{q\inv}  =
(\Ac_q)\inv$~\cite{GNS,selfsimilar}. Figure~\ref{aleshininverse}
shows the Moore diagram for the inverse of Aleshin's automaton from
Figure~\ref{aleshindiagram}.

It is well known that if $\mathcal A_q$ and $\mathcal B_s$ are
transformations computed by finite state initial automata, then the
composition $\mathcal B_s\mathcal A_q$ can be computed by a finite
state initial automaton~\cite{Eilenberg,GNS,selfsimilar}. If
$\mathcal A$ is an automaton over an alphabet $A$, we write
$\Sc(\mathcal A)$ for the semigroup of transformations of $A^*$ (or
equivalently $A^{\omega}$) generated by the functions $\Ac_q$ with
$q\in Q$. If $\mathcal A$ is invertible, we write $\GG (\Ac)$ for
the group of transformations generated by the initial automata
associated to the states of $\Ac_q$.  Groups generated by finite
state automata, also called \emph{automaton groups}, are a very
important special case of the general notion of a
\emph{self-similar} group~\cite{selfsimilar}.  A self-similar group
is what one gets by allowing infinite state automata in the
definition of an automaton group.  If $\Ac$ is an invertible
automaton, we write $\Ac^{\pm}$ for the automaton whose Moore
diagram is the disjoint union $\Ac$ and $\Ac^{\inv}$.  It is easy to
see that $\Ac^{\pm}$ is an invertible automaton, which is its own
inverse and that $\Sc(\Ac^{\pm})=\GG(\Ac)$.

If we denote by $T$ the rooted Cayley tree of $A^*$ with root vertex
the empty string, then $\GG(\Ac)$ acts on the left of $T$ by tree
automorphisms of $T$ \cite{branch,GNS,selfsimilar} via the action
\eqref{action}.  The induced action on the boundary $\partial T$
(the space of infinite directed paths from the root) is just the
action \eqref{actiononcantor} of $\GG(\Ac)$ on $A^{\omega}$.

The automorphism group $\mathrm{Aut}(T)$ is the iterated
(permutational) wreath product of countably many copies of the left
permutation group $(S_{A},A)$
\cite{Bass,branch,GNS,selfsimilar,Rhodestrees}, where $S_A$ denotes
the symmetric group on $A$.  In this paper, our notation will be
such that the wreath product of left permutation groups has a
natural projection to its leftmost factor; dual notation is used for
right permutation groups. For a group $\Ga=\GG (\Ac)$ generated by
an automaton over $A$, one has an embedding
\begin{equation}\label{wreathembedding1}
(\Ga,A^{\omega}) \hookrightarrow (S_{|A|},A)\wr (\Ga,A^{\omega})
\end{equation}
 where the
map sends $\Ac_q$ to the element with wreath product coordinates:
\begin{equation}\label{wreathcoords}
\Ac_q = \lambda_q(\Ac_{q_{a_1}},\ldots,\Ac_{q_{a_n}})
\end{equation}
where $A=\{a_1,\ldots,a_n\}$ and $\lambda_q(a) = \lambda(q,a)$.  See
\cite{branch,GNS,reset,selfsimilar} for more details. As an example,
if $\Ac$ is Aleshin's automaton from Figure~\ref{aleshindiagram},
then $a=(01)(c,b)$, $b= (01)(b,c)$ and $c=(a,a)$;  similarly, for
the inverse of Aleshin's automaton (Figure~\ref{aleshininverse}),
$a\inv = (01)(b\inv,c\inv)$, $b\inv = (01)(c\inv,b\inv)$ and $c\inv
= (a\inv,a\inv)$.  Notice that data \eqref{wreathcoords}, as $q$
varies over all states, completely encodes the automaton $\Ac$ and
so sometimes we shall give an automaton via the wreath product
coordinates.

We also shall need the notion of sections and the minimal automaton.
Let $a\in A$.
Then there is a homeomorphism $L_a:A^{\omega}\to A^{\omega}$ given
by $L_a(w)= aw$; metrically this is a contraction by a factor of
$1/|A|$.  If $f:A^{\omega}\to A^{\omega}$ is any metric
contraction, then, by definition of the metric, there is a function
$A\to A$, which we also denote by $f$ (abusing notation), so that,
for any word $w\in A^{\omega}$, $f(aw) = f(a)w'$, some $w'\in
A^{\omega}$. Define, for $a\in A$, $f_a:A^{\omega}\to A^{\omega}$ by
$f_a = L_{f(a)}\inv fL_a$. Then $f_a$ is also a metric contraction,
which is an isometry whenever $f$ is one.
The map $f_a$ is called the \emph{section} of
$f$ at $a$.  It is straightforward to verify that $f(aw) =
f(a)f_a(w)$. Notice that if $\mathcal A$ is a finite state
automaton, then $(\Ac_q)_a = \Ac_{q_a}$ --- that is the notation $q_a$
is unambiguous if we identify functions and states. One can also
define inductively the section $f_w$ for any word $w\in A^*$; of
course $f_{\varepsilon} =f$, where $\varepsilon$ is the empty string.
We remark that for any word $w\in A^*$ of length $n$, there is a
unique word of length $n$, denoted $f(w)$, so that
$f(wA^{\omega})\subseteq f(w)A^{\omega}$; this coincides with our
previous definition if $n=1$.  With this definition, one
has $f(wu) = f(w)f_w(u)$ for any word $u\in A^{\omega}$.

If $f:A^{\omega}\to A^{\omega}$ is a metric contraction, then the
minimal automaton of $f$ is the possibly infinite automaton $\Ac(f)$
with state set $\{f_w\mid w\in A^*\}$ (of course for different
$u,w\in A^*$, it may be the case that $f_u=f_w$).  The transitions
are given by $\delta(f_u,a) = f_{ua}$ and the output by
$\lambda(f_u,a)= f_u(a)$. It is easy to see that $\Ac(f)_{f_u} =
f_u$ and in particular $\Ac(f)_{f_{\varepsilon}} = f$.  One can
prove that $f$ is computed by a finite state automaton if and only
if $\Ac(f)$ is finite and that $\Ac(f)$ is the unique automaton with
minimal number of states computing $f$~\cite{Eilenberg,GNS}.
Moreover, it is well known that
$\Ac(f)$ is polynomial time computable from any automaton computing
$f$.  Notice
that if $f$ is invertible, then $\Ac(f)$ must be invertible and
$\Ac(f)\inv = \Ac(f\inv)$.

The following formula will be useful later.  Let
$f_n,\cdots,f_1:A^{\omega}\to A^{\omega}$ be metric contractions.
Then one easily checks that, for $w\in A^*$,
\begin{equation}\label{sectionofcomposition}
(f_n\cdots f_1)_w =  (f_n)_{f_{n-1}\cdots f_1(w)}\cdots
(f_2)_{f_1(w)}(f_1)_w
\end{equation}

Suppose that $\Ac=(Q,A^{\pm},\delta,\lambda)$ is an automaton and
$\sigma\in S_A$.  Then we need the following straightforward
observation on how to construct an automaton computing
$\sigma^*\Ac_q$ for any $q\in Q$.  The proof is left to the reader.

\begin{proposition}\label{composewithauto}
Let $\Ac=(Q,A,\delta,\lambda)$ be an automaton and $\sigma\in S_A$.
Define $\sigma[\Ac] = (Q,A,\delta,\sigma\lambda)$. Then
$\sigma[\Ac]_q =\sigma^*\Ac_q$.
\end{proposition}



In fact the proposition can be obtained as a special case of the
general construction for composing
automata~\cite{Eilenberg,GNS,selfsimilar}, using that automorphisms
are precisely the functions computed by invertible automata with a
single state.  If $A=Q^{\pm}$ and $\sigma\in S_A$, then in
particular we can apply the above proposition to understand
$\sigma^\pm\Ac_q$.

 Another notion that we shall
need is that of the transition monoid of an automaton.  If $\mathcal
A=(Q,A,\delta,\lambda)$ is an automaton, then the \emph{transition
monoid} of $\mathcal A$, denoted $M(\mathcal A)$, is the finite
monoid of all transformations of $Q$ of the form $q\mapsto q_w$ with
$w\in A^*$.  It is easy to check that this is indeed a finite
monoid.   It is well known~\cite{Eilenberg} that the transition
monoid $M(\Ac(f))$ of the minimal automaton $\mathcal A(f)$ is a
quotient of the transition monoid of any other initial automaton
computing $f$.  In other words, it is an algebraic invariant of $f$.
In this paper we shall be particularly interested in the case that
$M(\mathcal A)$ and $M(\mathcal A\inv)$ are groups.

A key fact about transition monoids is that if $f$ and $g$ are two
finite state transformations of $A^{\omega}$, then the transition
monoid $M(\mathcal A(fg))$ is a quotient of a submonoid of the
wreath product of right transformation monoids $M(\mathcal A(f))\wr
M(\mathcal A(g))$~\cite{Eilenberg,Arbib}.  Thus if $\mathcal C$ is a
class of finite monoids closed under taking wreath product,
submonoids and quotient monoids --- for instance the class of finite
groups --- then the set of finite state transformations with
transition monoid in $\mathcal C$ is a semigroup of metric
contractions.  We shall call the group of units of this semigroup
the \emph{group of $\mathcal C$-isometries of $A^{\omega}$}.  The
case where $\mathcal C$ is the collection of finite groups, gives
rise to the group of so-called bireversible
automata~\cite{glasner,MNS,free1,free2}, as we shall see below.

\subsection{The dual automaton}
The way we have defined automata, they scan their input from left to
right, outputting each time they process a letter.  Clearly one can
also define a \emph{right-to-left scanning automaton} in a dual
manner: so one has a $4$-tuple $(Q,A,\delta,\lambda)$ as before, but
now $\delta:A\times Q\to Q$ and $\lambda: A\times Q\to A$.  One then
sets $\delta(a,q) = {}_aq$ and $\lambda (a,q) = aq$ and extends to
$A^*$ via:

\begin{gather} {}_{ua}q = {}_u({}_aq) \\\label{dualaction} (ua)q =
u{}_aq(aq)
\end{gather}

If $\Ac$ is an invertible right-to-left scanning automaton, then
$\Ga=\GG(\Ac)$ is a group of permutations of $A^*$ (or isometries of
the space of left infinite words ${}^{\omega}\!{A}$, or
automorphisms of the left Cayley tree of $A^*$) acting on the right.
It preserves length and common suffixes, that is $(uw)\Ac_q =
u\Ac_{{}_wq}w\Ac_q$ for $u,w\in A^*$.

 There
is an embedding of $\Ga$ into
the wreath product of right permutation groups
\begin{equation}\label{wreathembeddingdual} (A^{\omega},\Ga)
\hookrightarrow ({}^{\omega}\!{A},\Ga)\wr (A,S_{|A|})
\end{equation}
 where the
map sends $\Ac_q$ to the element with wreath product coordinates:
\begin{equation}\label{wreathcoordsdual}
\Ac_q = (\Ac_{{}_{a_1}q},\ldots,\Ac_{{}_{a_n}q})\lambda_q
\end{equation}
where $A=\{a_1,\ldots,a_n\}$ and $a\lambda_q = \lambda(a,q)$.

Now if $\mathcal A=(Q,A,\delta,\lambda)$ is an automaton, then there
is an associated right-to-left scanning automaton, called the \emph{dual
automaton} of $\mathcal A$, which is denoted $\widehat{\Ac}$ and is
given by $\widehat{\Ac}=(A,Q,\lambda,\delta)$~\cite{GNS}.  So one
obtains $\widehat{\Ac}$ from $\mathcal A$ by switching the
states and the alphabet, switching the input and the output
functions and switching the way we scan words as the automaton
operates.  So one has for $q\in Q$ and $a\in A$
\begin{equation*}
\begin{split}
{}_qa &= q(a)\\ qa &= q_a
\end{split}
\end{equation*}

One then has, inductively, for $q_n,\ldots,q_1\in Q$ and $a\in A$,
\begin{equation}\label{dualize2}
\begin{split}
{}_{q_n\cdots q_2q_1}a &= q_nq_{n-1}\cdots q_1(a) = \Ac_{q_n}\cdots \Ac_{q_1}(a)\\
 q_n\cdots q_2q_1a &= (q_n)_{q_{n-1}\cdots
q_1(a)}\cdots (q_2)_{q_1(a)}(q_1)_a\\
&= (q_n)_{\Ac_{q_{n-1}}\cdots \Ac_{q_1}(a)}\cdots
(q_2)_{\Ac_{q_1}(a)}(q_1)_a
\end{split}
\end{equation}

One can interpret the first line of \eqref{dualize2} as saying that
in state $a$, on input $q_n\cdots q_1$, $\widehat{\Ac}$ goes to the
state which is the output of $\Ac_{q_n}\cdots \Ac_{q_1}$ on the
letter $a$. In other words, the transition diagram of $\wh{\Ac}$ is
the Schreier graph of the action of $\Sc(\Ac)$ on the first level of
the tree.  The second/third lines of \eqref{dualize2}, in light of
\eqref{sectionofcomposition}, says that the output of state $a$ with
input $q_n\cdots q_1$ is a word in $Q^*$ representing the section of
$\Ac_{q_n}\cdots \Ac_{q_1}$ at $a$.  The following result, proved
in~\cite{free1} in a slightly different language, is then immediate
by induction.

\begin{theorem}\label{dualprop}
Let $\mathcal A=(Q,A,\delta,\lambda)$ be a finite state automaton
with dual $\widehat{\Ac}$.  Then if $a_1,\ldots, a_m\in A$
and $q_n,\ldots,q_1\in Q$, then
\begin{align*}
{}_{q_n\cdots q_2q_1}\widehat{\Ac}_{a_1}\cdots \widehat{\Ac}_{a_m} & = \Ac_{q_n}\cdots
\Ac_{q_1}(a_1\cdots a_m)\\
(q_n\cdots q_2q_1)\widehat{\Ac}_{a_1}\cdots \widehat{\Ac}_{a_m}&=  (q_n)_{\Ac_{q_{n-1}}\cdots \Ac_{q_1}(a_1\cdots a_m)}\cdots
(q_2)_{\Ac_{q_1}(a_1\cdots a_m)}(q_1)_{a_1\cdots a_m}
\end{align*}
In particular, by \eqref{sectionofcomposition}, $(q_n\cdots
q_2q_1)\widehat{\Ac}_{a_1}\cdots \widehat{\Ac}_{a_m}$ is a word in
$Q^*$ representing in $\Sc(\Ac)$ the section $(\Ac_{q_n}\cdots
\Ac_{q_1})_{a_1\cdots a_m}$ of $q_n\cdots q_1$ at $a_1\cdots a_m$.
\end{theorem}

\begin{corollary}\label{orbitcondition}
Let $\Ac$ be an automaton with state set $Q$.  Suppose that $u,v\in
Q^*$ and $v$ is in the orbit of $u$ under $\Sc(\widehat{\Ac})$, that
is $v=us$ for some $s\in \Sc(\widehat{\Ac})$.  Then $v$ represents a
section of $u$.
\end{corollary}

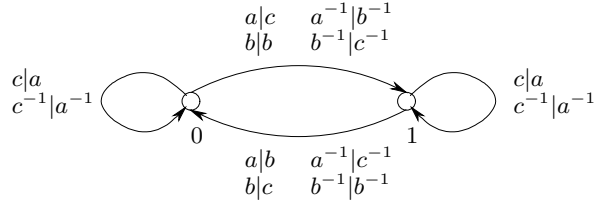
\begin{figure}[htbp]
\begin{center}
\setlength{\unitlength}{0.00062500in}
\begingroup\makeatletter\ifx\SetFigFont\undefined%
\gdef\SetFigFont#1#2#3#4#5{%
  \reset@font\fontsize{#1}{#2pt}%
  \fontfamily{#3}\fontseries{#4}\fontshape{#5}%
  \selectfont}%
\fi\endgroup%
\begin{picture}(4892,1575)(0,-10)
\put(2400.000,-787.500){\arc{3825.000}{4.2224}{5.2023}}
\blacken\path(3179.145,926.343)(3300.000,900.000)(3205.760,980.117)(3224.717,937.261)(3179.145,926.343)
\put(2400.000,2437.000){\arc{3824.118}{1.0807}{2.0609}}
\blacken\path(1620.853,723.645)(1500.000,750.000)(1594.233,669.874)(1575.280,712.732)(1620.853,723.645)
\put(1500,825){\ellipse{150}{150}}
\put(3303,823){\ellipse{150}{150}}
\path(1464,872)(1463,873)(1459,877)
    (1454,882)(1445,890)(1433,901)
    (1419,914)(1402,929)(1383,944)
    (1363,961)(1341,977)(1318,993)
    (1294,1008)(1269,1023)(1242,1035)
    (1213,1046)(1183,1055)(1150,1062)
    (1116,1065)(1080,1065)(1048,1061)
    (1018,1055)(990,1047)(965,1039)
    (943,1030)(923,1022)(906,1013)
    (891,1005)(878,998)(866,990)
    (854,982)(843,975)(832,966)
    (822,957)(810,947)(799,935)
    (788,921)(777,906)(766,888)
    (758,869)(752,847)(750,825)
    (753,803)(760,781)(770,761)
    (782,743)(794,726)(806,712)
    (818,699)(830,688)(842,678)
    (853,668)(865,659)(877,651)
    (890,642)(904,634)(919,625)
    (937,615)(957,606)(979,596)
    (1005,587)(1033,579)(1063,572)
    (1095,568)(1130,568)(1164,571)
    (1195,578)(1225,588)(1252,600)
    (1277,613)(1301,628)(1323,644)
    (1344,662)(1364,679)(1383,697)
    (1400,714)(1415,729)(1428,743)
    (1439,755)(1456,774)
\blacken\path(1398.342,664.567)(1456.000,774.000)(1353.627,704.575)(1399.989,711.400)(1398.342,664.567)
\path(3336,872)(3337,873)(3341,877)
    (3346,882)(3355,890)(3367,901)
    (3381,914)(3398,929)(3417,944)
    (3437,961)(3459,977)(3482,993)
    (3506,1008)(3531,1023)(3558,1035)
    (3587,1046)(3617,1055)(3650,1062)
    (3684,1065)(3720,1065)(3752,1061)
    (3782,1055)(3810,1047)(3835,1039)
    (3857,1030)(3877,1022)(3894,1013)
    (3909,1005)(3922,998)(3934,990)
    (3946,982)(3957,975)(3968,966)
    (3978,957)(3990,947)(4001,935)
    (4012,921)(4023,906)(4034,888)
    (4042,869)(4048,847)(4050,825)
    (4047,803)(4040,781)(4030,761)
    (4018,743)(4006,726)(3994,712)
    (3982,699)(3970,688)(3958,678)
    (3947,668)(3935,659)(3923,651)
    (3910,642)(3896,634)(3881,625)
    (3863,615)(3843,606)(3821,596)
    (3795,587)(3767,579)(3737,572)
    (3705,568)(3670,568)(3636,571)
    (3605,578)(3575,588)(3548,600)
    (3523,613)(3499,628)(3477,644)
    (3456,662)(3436,679)(3417,697)
    (3400,714)(3385,729)(3372,743)
    (3361,755)(3344,774)
\blacken\path(3446.373,704.575)(3344.000,774.000)(3401.658,664.567)(3400.011,711.400)(3446.373,704.575)
\put(3300,475){\makebox(0,0)[lb]{{\SetFigFont{9}{10.8}{\rmdefault}{\mddefault}{\updefault}$1$}}}
\put(1500,475){\makebox(0,0)[lb]{{\SetFigFont{9}{10.8}{\rmdefault}{\mddefault}{\updefault}$0$}}}
\put(1950,1425){\makebox(0,0)[lb]{{\SetFigFont{9}{10.8}{\rmdefault}{\mddefault}{\updefault}$a|c$}}}
\put(2500,1425){\makebox(0,0)[lb]{{\SetFigFont{9}{10.8}{\rmdefault}{\mddefault}{\updefault}$a\inv|b\inv$}}}
\put(1950,225){\makebox(0,0)[lb]{{\SetFigFont{9}{10.8}{\rmdefault}{\mddefault}{\updefault}$a|b$}}}
\put(2500,225){\makebox(0,0)[lb]{{\SetFigFont{9}{10.8}{\rmdefault}{\mddefault}{\updefault}$a\inv|c\inv$}}}
\put(4200,900){\makebox(0,0)[lb]{{\SetFigFont{9}{10.8}{\rmdefault}{\mddefault}{\updefault}$c|a$}}}
\put(4200,675){\makebox(0,0)[lb]{{\SetFigFont{9}{10.8}{\rmdefault}{\mddefault}{\updefault}$c\inv|a\inv$}}}
\put(1950,1200){\makebox(0,0)[lb]{{\SetFigFont{9}{10.8}{\rmdefault}{\mddefault}{\updefault}$b|b$}}}
\put(2500,1200){\makebox(0,0)[lb]{{\SetFigFont{9}{10.8}{\rmdefault}{\mddefault}{\updefault}$b\inv|c\inv$}}}
\put(1950,0){\makebox(0,0)[lb]{{\SetFigFont{9}{10.8}{\rmdefault}{\mddefault}{\updefault}$b|c$}}}
\put(2500,0){\makebox(0,0)[lb]{{\SetFigFont{9}{10.8}{\rmdefault}{\mddefault}{\updefault}$b\inv|b\inv$}}}
\put(0,675){\makebox(0,0)[lb]{{\SetFigFont{9}{10.8}{\rmdefault}{\mddefault}{\updefault}$c\inv|a\inv$}}}
\put(0,900){\makebox(0,0)[lb]{{\SetFigFont{9}{10.8}{\rmdefault}{\mddefault}{\updefault}$c|a$}}}
\end{picture}
\end{center}
\caption{The dual of $\mathsf A^{\pm}$\label{aleshindual}}
\end{figure}

Let $\mathsf A$ be the Aleshin automaton from Figure~\ref{aleshindiagram}.
Then the dual automaton for $\mathsf A^{\pm}$ is given by
Figure~\ref{aleshindual}.

\subsection{Reversible and bireversible automata}
The dual automaton is most useful when it is invertible.  We discuss
this situation here.  A finite state automaton $\mathcal
A=(Q,A,\delta,\lambda)$ is called \emph{reversible} if, for all
$a\in A$, $q,q'\in Q$, one has that $q_a=q'_a$ implies $q=q'$.  That
is the map $q\to q_a$ is a permutation for all $a\in A$.  This is
equivalent to asking that the transition monoid $M(\Ac)$ be a group.
In particular, a function $f:A^{\omega}\to A^{\omega}$ can be
computed by a reversible automaton if and only if the minimal
automaton $\Ac(f)$ is reversible.  We then say that the metric
contraction $f$ is \emph{reversible}.   An invertible automaton is
called \emph{bireversible} if both it and its inverse are
reversible.  Notice that if $\Ac$ is bireversible, then $\Ac^{\pm}$
is also bireversible. An isometry $f:A^{\omega}\to A^{\omega}$ is
called bireversible if $f$ and $f^{-1}$ are reversible or,
equivalently, the minimal automaton $\Ac(f)$ is bireversible.

As mentioned earlier, since the class of finite groups is closed
under wreath product, submonoids and quotients, the collections of
all reversible metric contractions of $A^{\omega}$ forms a
semigroup, denoted $\mathsf {Rev}(A)$, called the semigroup of
reversible automata.  The group of units of $\mathsf {Rev}(A)$ is
then denoted $\mathsf{BiRev}(A)$ and its elements are precisely the
bireversible isometries.   One calls $\mathsf{BiRev}(A)$ the group
of bireversible automata. See~\cite{glasner,MNS} for relations
between $\mathsf{BiRev}(A)$ and commensurators.

If $\Ac$ is a bireversible automaton, then $\GG(\Ac)$ is a finitely
generated subgroup of $\mathsf{BiRev}(A)$.  It turns out that such
groups are particularly apt for analysis via their dual automata and
they have some remarkable properties. First of all, the
reversibility of $\Ac$ and $\Ac\inv$ is equivalent to asking that,
for each $a\in A$, the maps $q\mapsto q_a$ and $q\inv \mapsto
q\inv_a$ are invertible. But this is exactly the same as asking that
$\widehat{\Ac}$ and $\widehat{\Ac\inv}$ be invertible.  This proves
the following well-known result~\cite{selfsimilar,free1}.

\begin{proposition}\label{bireversibleintermsofdual}
A finite automaton $\mathcal A$ is bireversible if and only if
$\widehat{\Ac}$ and $\widehat{\Ac\inv}$ are invertible.
\end{proposition}

We also need the following trivial observation:  any subsemigroup of
a finite group is a subgroup since the inverse of an element is a
positive power.  In particular, any semigroup of permutations of a
finite set is a group.  If $\Ac$ is an invertible automaton acting
on $A^*$, then $\Sc(\Ac)$ acts by permutations on the finite set
$A^n$ and hence acts as a finite group of permutations of $A^n$.
Thus $\Sc(\Ac)$ and $\GG(\Ac)$ have the same orbits on $A^n$.
Summarizing, we obtain the following~\cite{free1}.

\begin{lemma}\label{semigrouporbits}
Let $\mathcal A=(Q,A,\delta,\lambda)$ be an invertible automaton.
Then the orbits of $\Sc(\mathcal A)$ and of $\GG(\mathcal A)$ on $A^*$
are the same.
\end{lemma}

\begin{corollary}\label{bireversiblesectioncor}
Let $\mathcal A=(Q,A,\delta,\lambda)$ be a bireverisble automaton.
Let $u,v\in Q^*$ and suppose that $ug=v$ for some $g\in
\GG(\widehat{\Ac})$.  Then $v$ represents a section of $u$ and $u$
represents a section of $v$.
\end{corollary}
\begin{proof}
Since $u=vg\inv$, it suffices to show that $v$ represents a section of
$u$.  But this follows from Lemma~\ref{semigrouporbits} and
Corollary~\ref{dualprop}.
\end{proof}

In particular, we have the following criterion for triviality in a
bireversible automaton group:

\begin{corollary}[Triviality criterion]\label{triviality}
Let $\Ac=(Q,A,\delta,\lambda)$ be a bireversible automaton and let
$w\in (Q^\pm)^*$.  Then the following are equivalent:
\begin{enumerate}
\item $w$ represents the trivial element of $\GG(\Ac)$
\item all sections of $w$ are trivial
\item $w$ has a section that is trivial
\item $w$ is in the orbit under $\GG(\widehat{\Ac^\pm})$ of an
  element $u\in (Q^\pm)^*$ that is trivial in $\GG(\Ac)$
\item the whole orbit of $w$ under $\GG(\widehat{\Ac^\pm})$
  represents the trivial element of $\GG(\Ac)$.
\end{enumerate}
\end{corollary}
\begin{proof}
Clearly (1) implies (2) since all sections of the identity
transformation are trivial.  The implication $(2)\implies (3)$ is
trivial.  For $(3)\implies (4)$,  suppose that the section of $w$ at
$a_1\cdots a_m\in A^*$ is trivial.  Then, by Theorem~\ref{dualprop},
$u= w\widehat{\Ac^\pm}_{a_1}\cdots \widehat{\Ac^\pm}_{a_m}$
represents the section of $w$ at $a_1\cdots a_m$ and hence is
trivial in $\GG(\Ac)$, so (4) holds.  To see that (4) implies (1),
we have by Corollary~\ref{bireversiblesectioncor} that $w$ is a
section of the transformation represented by $u$.  But since $u$
represents the identity transformation and every section of the
identity is the identity, it follows that $w$ represents the trivial
element of $\GG(\Ac)$.  Clearly (5) implies (4).  On the other hand,
(2) implies (5) by Corollary~\ref{bireversiblesectioncor}.
\end{proof}

The implication (4) implies (1) was established in~\cite{free1} and is
the key tool to proving freeness.  Indeed, there is the following
criterion for freeness of a bireversible automaton group that is
immediate from Corollary~\ref{triviality}.

\begin{corollary}[Freeness criterion]\label{freeness}
Let $\Ac=(Q,A,\delta,\lambda)$ be a bireversible automaton.  Then
$\GG(\Ac)$ is a free group, freely generated by the transformations
$\Ac_q$, $q\in Q$, if and only for each freely irreducible word
$w\in (Q^\pm)^*$, there is a word $u$ in the orbit of $w$ under
$\GG(\widehat{\Ac^\pm})$ such that $u$ represents a non-trivial
element of $\GG(\Ac)$.
\end{corollary}

This leads to the following basic strategy for proving that the group
generated by a bireversible automaton is free.

\begin{definition}[Pattern]
A word in the alphabet $\{\ast,\ast\inv\}$ is called a \emph{pattern}.   A
word $w\in (Q^\pm)^*$ is said to \emph{follow} the
pattern $u$ if the image of $w$ under the map sending $Q$ to $\ast$
and $Q\inv$ to $\ast\inv$ is $u$.
\end{definition}

The basic strategy is then to show that all freely irreducible words following
any given pattern are in the same orbit of $\GG(\widehat{\Ac^{\pm}})$ and then
to show that each pattern is followed by some element
that does not act trivially on words of length $1$.
This strategy was successfully
used by the second and third authors to prove that the Aleshin
automata and its relatives are free~\cite{free1,free2}.  Glasner and
Mozes used a similar criterion in their construction of free groups
generated by automata~\cite{glasner}.

Let us point out another property of bireversible automata that is
implied by Corollary~\ref{triviality}.  We view $A^{\omega}$ as a
measure space by taking the product of the uniform
measure on $A$.  For a transformation $f:A^{\omega}\to A^{\omega}$, we
denote by $\mathrm{Fix}(f)$ the set of fixed points of $f$; it is a
closed subspace of $A^{\omega}$.  The reader is referred
to~\cite{GZihara,kss} for the definition of the Kesten and the Kesten-von
Neumann-Serre spectral measures.

\begin{corollary}
Let $\Ac=(Q,A,\delta,\lambda)$ be a bireversible automaton and let
$1\neq g\in \GG(\Ac)$.  The $\mathrm{Fix}(g)$ is nowhere dense and has
measure zero.   Hence the Kesten and the Kesten-von
Neumann-Serre spectral measures with respect to any generating set for
$\GG(\Ac)$ coincide.
\end{corollary}
\begin{proof}
Suppose that $\mathrm{Fix}(g)$ is not nowhere
dense.  Then it contains a cylinder set $wA^{\omega}$ with $w\in
A^*$.  It follows that the section $g_w$ is trivial.  But then
Corollary~\ref{triviality} shows that $g$ is trivial.

It is shown
in~\cite{kss} that, for any transformation $f:A^{\omega}\to A^{\omega}$
computed by a finite state automaton, $\mathrm{Fix}(f)$ is nowhere
dense if and only if it has measure zero, proving the second
statement.

The condition  $\mathrm{Fix}(g)$ has measure zero for all
non-trivial elements $g\in \GG(\Ac)$ was shown in~\cite{kss} to
imply that the Kesten and the Kesten-von Neumann-Serre spectral
measures with respect to any generating set for $\GG(\Ac)$ coincide.
\end{proof}

This means that for bireversible automata, one can attempt to compute
the Kesten spectral measure along the lines of~\cite{GZlamp,kss}.
Having completed the general theory we need for this paper, we now turn
to the series of finite state automata in question.

\section{The dramatis person\ae}
In this section we shall introduce the various automata used
throughout this paper.  We shall try to remind the reader explicitly
when an automaton is right-to-left scanning, but it should be kept
that we always take this to be the case for dual automata.

\subsection{The Aleshin automaton and its relatives}
We record here several results from~\cite{free1} about Aleshin's
automaton that we shall require
in the sequel.  Let
$\mathsf A$ be Aleshin's automaton from Figure~\ref{aleshindiagram}
and let $\mathsf D=\wh{\mathsf A^{\pm}}$ be the right-to-left scanning automaton from
Figure~\ref{aleshindual}.

\begin{theorem}[\cite{free1}]\label{mainaleshinresult}
If $p\in \{\ast,\ast\inv\}^*$ is a pattern, then $\GG(\mathsf D)$
acts transitively on the set of freely irreducible words following
$p$.  In particular, $\GG(\mathsf D)$ acts transitively on
$\{a,b,c\}^*$.
\end{theorem}

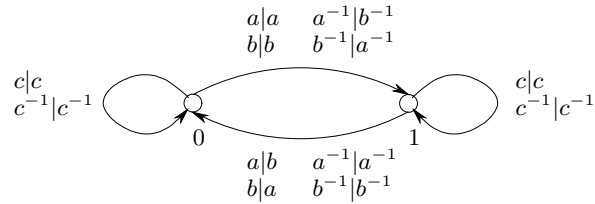
\begin{figure}[htbp]
\begin{center}
\setlength{\unitlength}{0.00062500in}
\begingroup\makeatletter\ifx\SetFigFont\undefined%
\gdef\SetFigFont#1#2#3#4#5{%
  \reset@font\fontsize{#1}{#2pt}%
  \fontfamily{#3}\fontseries{#4}\fontshape{#5}%
  \selectfont}%
\fi\endgroup%
\begin{picture}(4892,1575)(0,-10)
\put(2400.000,-787.500){\arc{3825.000}{4.2224}{5.2023}}
\blacken\path(3179.145,926.343)(3300.000,900.000)(3205.760,980.117)(3224.717,937.261)(3179.145,926.343)
\put(2400.000,2437.000){\arc{3824.118}{1.0807}{2.0609}}
\blacken\path(1620.853,723.645)(1500.000,750.000)(1594.233,669.874)(1575.280,712.732)(1620.853,723.645)
\put(1500,825){\ellipse{150}{150}}
\put(3303,823){\ellipse{150}{150}}
\path(1464,872)(1463,873)(1459,877)
    (1454,882)(1445,890)(1433,901)
    (1419,914)(1402,929)(1383,944)
    (1363,961)(1341,977)(1318,993)
    (1294,1008)(1269,1023)(1242,1035)
    (1213,1046)(1183,1055)(1150,1062)
    (1116,1065)(1080,1065)(1048,1061)
    (1018,1055)(990,1047)(965,1039)
    (943,1030)(923,1022)(906,1013)
    (891,1005)(878,998)(866,990)
    (854,982)(843,975)(832,966)
    (822,957)(810,947)(799,935)
    (788,921)(777,906)(766,888)
    (758,869)(752,847)(750,825)
    (753,803)(760,781)(770,761)
    (782,743)(794,726)(806,712)
    (818,699)(830,688)(842,678)
    (853,668)(865,659)(877,651)
    (890,642)(904,634)(919,625)
    (937,615)(957,606)(979,596)
    (1005,587)(1033,579)(1063,572)
    (1095,568)(1130,568)(1164,571)
    (1195,578)(1225,588)(1252,600)
    (1277,613)(1301,628)(1323,644)
    (1344,662)(1364,679)(1383,697)
    (1400,714)(1415,729)(1428,743)
    (1439,755)(1456,774)
\blacken\path(1398.342,664.567)(1456.000,774.000)(1353.627,704.575)(1399.989,711.400)(1398.342,664.567)
\path(3336,872)(3337,873)(3341,877)
    (3346,882)(3355,890)(3367,901)
    (3381,914)(3398,929)(3417,944)
    (3437,961)(3459,977)(3482,993)
    (3506,1008)(3531,1023)(3558,1035)
    (3587,1046)(3617,1055)(3650,1062)
    (3684,1065)(3720,1065)(3752,1061)
    (3782,1055)(3810,1047)(3835,1039)
    (3857,1030)(3877,1022)(3894,1013)
    (3909,1005)(3922,998)(3934,990)
    (3946,982)(3957,975)(3968,966)
    (3978,957)(3990,947)(4001,935)
    (4012,921)(4023,906)(4034,888)
    (4042,869)(4048,847)(4050,825)
    (4047,803)(4040,781)(4030,761)
    (4018,743)(4006,726)(3994,712)
    (3982,699)(3970,688)(3958,678)
    (3947,668)(3935,659)(3923,651)
    (3910,642)(3896,634)(3881,625)
    (3863,615)(3843,606)(3821,596)
    (3795,587)(3767,579)(3737,572)
    (3705,568)(3670,568)(3636,571)
    (3605,578)(3575,588)(3548,600)
    (3523,613)(3499,628)(3477,644)
    (3456,662)(3436,679)(3417,697)
    (3400,714)(3385,729)(3372,743)
    (3361,755)(3344,774)
\blacken\path(3446.373,704.575)(3344.000,774.000)(3401.658,664.567)(3400.011,711.400)(3446.373,704.575)
\put(3300,475){\makebox(0,0)[lb]{{\SetFigFont{9}{10.8}{\rmdefault}{\mddefault}{\updefault}$1$}}}
\put(1500,475){\makebox(0,0)[lb]{{\SetFigFont{9}{10.8}{\rmdefault}{\mddefault}{\updefault}$0$}}}
\put(1950,1425){\makebox(0,0)[lb]{{\SetFigFont{9}{10.8}{\rmdefault}{\mddefault}{\updefault}$a|a$}}}
\put(2500,1425){\makebox(0,0)[lb]{{\SetFigFont{9}{10.8}{\rmdefault}{\mddefault}{\updefault}$a\inv|b\inv$}}}
\put(1950,225){\makebox(0,0)[lb]{{\SetFigFont{9}{10.8}{\rmdefault}{\mddefault}{\updefault}$a|b$}}}
\put(2500,225){\makebox(0,0)[lb]{{\SetFigFont{9}{10.8}{\rmdefault}{\mddefault}{\updefault}$a\inv|a\inv$}}}
\put(4200,900){\makebox(0,0)[lb]{{\SetFigFont{9}{10.8}{\rmdefault}{\mddefault}{\updefault}$c|c$}}}
\put(4200,675){\makebox(0,0)[lb]{{\SetFigFont{9}{10.8}{\rmdefault}{\mddefault}{\updefault}$c\inv|c\inv$}}}
\put(1950,1200){\makebox(0,0)[lb]{{\SetFigFont{9}{10.8}{\rmdefault}{\mddefault}{\updefault}$b|b$}}}
\put(2500,1200){\makebox(0,0)[lb]{{\SetFigFont{9}{10.8}{\rmdefault}{\mddefault}{\updefault}$b\inv|a\inv$}}}
\put(1950,0){\makebox(0,0)[lb]{{\SetFigFont{9}{10.8}{\rmdefault}{\mddefault}{\updefault}$b|a$}}}
\put(2500,0){\makebox(0,0)[lb]{{\SetFigFont{9}{10.8}{\rmdefault}{\mddefault}{\updefault}$b\inv|b\inv$}}}
\put(0,675){\makebox(0,0)[lb]{{\SetFigFont{9}{10.8}{\rmdefault}{\mddefault}{\updefault}$c\inv|c\inv$}}}
\put(0,900){\makebox(0,0)[lb]{{\SetFigFont{9}{10.8}{\rmdefault}{\mddefault}{\updefault}$c|c$}}}
\end{picture}
\end{center}
\caption{The right-to-left scanning automaton $\mathsf E$\label{aleshinE}}
\end{figure}

In the paper~\cite{free1}, the dual automaton is viewed as scanning
from left to right, instead of right to left as we do in this
paper.  But this does not affect the validity of
Theorem~\ref{mainaleshinresult} since the notions of pattern and
freely irreducible words are left-right dual, as is $\{a,b,c\}^*$.
Consider the right-to-left scanning automaton $\mathsf E$ whose Moore diagram is in
Figure~\ref{aleshinE}.

\begin{proposition}[\cite{free1}]\label{generatorsforaleshindual}
$\GG(\mathsf D)$ is generated by $\mathsf E_0$ and the automorphisms
$(ab)^{\pm}$, $(bc)^{\pm}$.
\end{proposition}

We remark that $\mathsf E_1=\mathsf E_0(ab)^\pm$~\cite{free1} and so
also belongs to $\GG(\mathsf D)$.

\subsection{A family of bireversible automata with an even number of states}
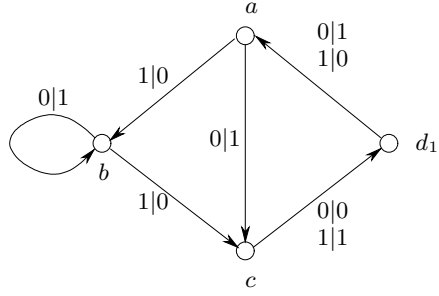
\begin{figure}[htbp]
\begin{center}
\setlength{\unitlength}{0.00062500in}
\begingroup\makeatletter\ifx\SetFigFont\undefined%
\gdef\SetFigFont#1#2#3#4#5{%
  \reset@font\fontsize{#1}{#2pt}%
  \fontfamily{#3}\fontseries{#4}\fontshape{#5}%
  \selectfont}%
\fi\endgroup%
\begin{picture}(3593,2395)(0,-10)
\put(784,1200){\ellipse{150}{150}}
\put(1984,2100){\ellipse{150}{150}}
\put(1984,300){\ellipse{150}{150}}
\put(3184,1200){\ellipse{150}{150}}
\path(856,1152)(1921,327)
\blacken\path(1807.762,376.771)(1921.000,327.000)(1844.506,424.204)(1854.594,378.441)(1807.762,376.771)
\path(1893,2075)(853,1227)
\blacken\path(927.044,1326.083)(853.000,1227.000)(964.960,1279.582)(918.102,1280.083)(927.044,1326.083)
\path(2051,330)(3116,1155)
\blacken\path(3039.506,1057.796)(3116.000,1155.000)(3002.762,1105.229)(3049.594,1103.559)(3039.506,1057.796)
\path(3115,1247)(2075,2095)
\blacken\path(2186.960,2042.418)(2075.000,2095.000)(2149.044,1995.917)(2140.102,2041.917)(2186.960,2042.418)
\path(1984,2025)(1984,375)
\blacken\path(1954.000,495.000)(1984.000,375.000)(2014.000,495.000)(1984.000,459.000)(1954.000,495.000)
\path(726,1246)(725,1247)(721,1251)
    (716,1256)(707,1264)(695,1275)
    (681,1288)(664,1303)(645,1318)
    (625,1335)(603,1351)(580,1367)
    (556,1382)(531,1397)(504,1409)
    (475,1420)(445,1429)(412,1436)
    (378,1439)(342,1439)(310,1435)
    (280,1429)(252,1421)(227,1413)
    (205,1404)(185,1396)(168,1387)
    (153,1379)(140,1372)(128,1364)
    (116,1356)(105,1349)(94,1340)
    (84,1331)(72,1321)(61,1309)
    (50,1295)(39,1280)(28,1262)
    (20,1243)(14,1221)(12,1199)
    (15,1177)(22,1155)(32,1135)
    (44,1117)(56,1100)(68,1086)
    (80,1073)(92,1062)(104,1052)
    (115,1042)(127,1033)(139,1025)
    (152,1016)(166,1008)(181,999)
    (199,989)(219,980)(241,970)
    (267,961)(295,953)(325,946)
    (357,942)(392,942)(426,945)
    (457,952)(487,962)(514,974)
    (539,987)(563,1002)(585,1018)
    (606,1036)(626,1053)(645,1071)
    (662,1088)(677,1103)(690,1117)
    (701,1129)(718,1148)
\blacken\path(660.342,1038.567)(718.000,1148.000)(615.627,1078.575)(661.989,1085.400)(660.342,1038.567)
\put(1984,0){\makebox(0,0)[lb]{{\SetFigFont{9}{10.8}{\rmdefault}{\mddefault}{\updefault}$c$}}}
\put(754,890){\makebox(0,0)[lb]{{\SetFigFont{9}{10.8}{\rmdefault}{\mddefault}{\updefault}$b$}}}
\put(1984,2290){\makebox(0,0)[lb]{{\SetFigFont{9}{10.8}{\rmdefault}{\mddefault}{\updefault}$a$}}}
\put(252,1467){\makebox(0,0)[lb]{{\SetFigFont{9}{10.8}{\rmdefault}{\mddefault}{\updefault}$0|1$}}}
\put(1082,1655){\makebox(0,0)[lb]{{\SetFigFont{9}{10.8}{\rmdefault}{\mddefault}{\updefault}$1|0$}}}
\put(1082,605){\makebox(0,0)[lb]{{\SetFigFont{9}{10.8}{\rmdefault}{\mddefault}{\updefault}$1|0$}}}
\put(2584,525){\makebox(0,0)[lb]{{\SetFigFont{9}{10.8}{\rmdefault}{\mddefault}{\updefault}$0|0$}}}
\put(2584,300){\makebox(0,0)[lb]{{\SetFigFont{9}{10.8}{\rmdefault}{\mddefault}{\updefault}$1|1$}}}
\put(2584,2025){\makebox(0,0)[lb]{{\SetFigFont{9}{10.8}{\rmdefault}{\mddefault}{\updefault}$0|1$}}}
\put(2584,1800){\makebox(0,0)[lb]{{\SetFigFont{9}{10.8}{\rmdefault}{\mddefault}{\updefault}$1|0$}}}
\put(1684,1125){\makebox(0,0)[lb]{{\SetFigFont{9}{10.8}{\rmdefault}{\mddefault}{\updefault}$0|1$}}}
\put(3409,1125){\makebox(0,0)[lb]{{\SetFigFont{9}{10.8}{\rmdefault}{\mddefault}{\updefault}$d_1$}}}
\end{picture}
\end{center}
\caption{A four-state automaton from $\mathfrak F$\label{fourstate}}
\end{figure}

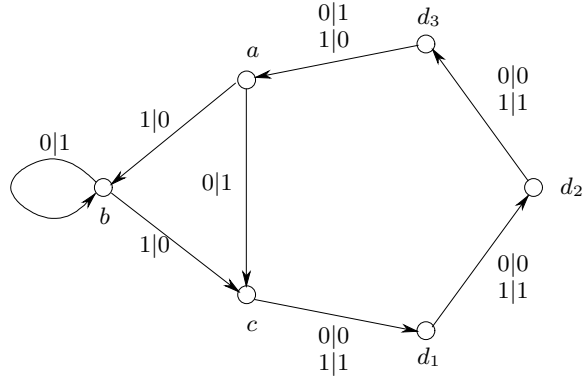
\begin{figure}[htbp]
\begin{center}
\setlength{\unitlength}{0.00062500in}
\begingroup\makeatletter\ifx\SetFigFont\undefined%
\gdef\SetFigFont#1#2#3#4#5{%
  \reset@font\fontsize{#1}{#2pt}%
  \fontfamily{#3}\fontseries{#4}\fontshape{#5}%
  \selectfont}%
\fi\endgroup%
\begin{picture}(4793,3075)(0,-10)
\put(2584,225){\makebox(0,0)[lb]{{\SetFigFont{9}{10.8}{\rmdefault}{\mddefault}{\updefault}$0|0$}}}
\put(2584,0){\makebox(0,0)[lb]{{\SetFigFont{9}{10.8}{\rmdefault}{\mddefault}{\updefault}$1|1$}}}
\put(4084,825){\makebox(0,0)[lb]{{\SetFigFont{9}{10.8}{\rmdefault}{\mddefault}{\updefault}$0|0$}}}
\put(4084,600){\makebox(0,0)[lb]{{\SetFigFont{9}{10.8}{\rmdefault}{\mddefault}{\updefault}$1|1$}}}
\put(4084,2400){\makebox(0,0)[lb]{{\SetFigFont{9}{10.8}{\rmdefault}{\mddefault}{\updefault}$0|0$}}}
\put(4084,2175){\makebox(0,0)[lb]{{\SetFigFont{9}{10.8}{\rmdefault}{\mddefault}{\updefault}$1|1$}}}
\put(784,1575){\ellipse{150}{150}}
\put(1984,2475){\ellipse{150}{150}}
\put(1984,675){\ellipse{150}{150}}
\put(3484,375){\ellipse{150}{150}}
\put(4384,1575){\ellipse{150}{150}}
\put(3484,2775){\ellipse{150}{150}}
\path(856,1527)(1921,702)
\blacken\path(1807.762,751.771)(1921.000,702.000)(1844.506,799.204)(1854.594,753.441)(1807.762,751.771)
\path(1893,2450)(853,1602)
\blacken\path(927.044,1701.083)(853.000,1602.000)(964.960,1654.582)(918.102,1655.083)(927.044,1701.083)
\path(3539,408)(4339,1508)
\blacken\path(4292.681,1393.307)(4339.000,1508.000)(4244.157,1428.597)(4289.594,1440.066)(4292.681,1393.307)
\path(4341,1631)(3541,2731)
\blacken\path(3635.843,2651.597)(3541.000,2731.000)(3587.319,2616.307)(3590.406,2663.066)(3635.843,2651.597)
\path(3421,2765)(2060,2499)
\blacken\path(2172.017,2551.461)(2060.000,2499.000)(2183.526,2492.575)(2142.440,2515.112)(2172.017,2551.461)
\path(2053,638)(3414,372)
\blacken\path(3290.474,365.575)(3414.000,372.000)(3301.983,424.461)(3331.560,388.112)(3290.474,365.575)
\path(1984,2400)(1984,750)
\blacken\path(1954.000,870.000)(1984.000,750.000)(2014.000,870.000)(1984.000,834.000)(1954.000,870.000)
\path(726,1621)(725,1622)(721,1626)
    (716,1631)(707,1639)(695,1650)
    (681,1663)(664,1678)(645,1693)
    (625,1710)(603,1726)(580,1742)
    (556,1757)(531,1772)(504,1784)
    (475,1795)(445,1804)(412,1811)
    (378,1814)(342,1814)(310,1810)
    (280,1804)(252,1796)(227,1788)
    (205,1779)(185,1771)(168,1762)
    (153,1754)(140,1747)(128,1739)
    (116,1731)(105,1724)(94,1715)
    (84,1706)(72,1696)(61,1684)
    (50,1670)(39,1655)(28,1637)
    (20,1618)(14,1596)(12,1574)
    (15,1552)(22,1530)(32,1510)
    (44,1492)(56,1475)(68,1461)
    (80,1448)(92,1437)(104,1427)
    (115,1417)(127,1408)(139,1400)
    (152,1391)(166,1383)(181,1374)
    (199,1364)(219,1355)(241,1345)
    (267,1336)(295,1328)(325,1321)
    (357,1317)(392,1317)(426,1320)
    (457,1327)(487,1337)(514,1349)
    (539,1362)(563,1377)(585,1393)
    (606,1411)(626,1428)(645,1446)
    (662,1463)(677,1478)(690,1492)
    (701,1504)(718,1523)
\blacken\path(660.342,1413.567)(718.000,1523.000)(615.627,1453.575)(661.989,1460.400)(660.342,1413.567)
\put(1984,375){\makebox(0,0)[lb]{{\SetFigFont{9}{10.8}{\rmdefault}{\mddefault}{\updefault}$c$}}}
\put(754,1265){\makebox(0,0)[lb]{{\SetFigFont{9}{10.8}{\rmdefault}{\mddefault}{\updefault}$b$}}}
\put(1984,2665){\makebox(0,0)[lb]{{\SetFigFont{9}{10.8}{\rmdefault}{\mddefault}{\updefault}$a$}}}
\put(252,1842){\makebox(0,0)[lb]{{\SetFigFont{9}{10.8}{\rmdefault}{\mddefault}{\updefault}$0|1$}}}
\put(1082,2030){\makebox(0,0)[lb]{{\SetFigFont{9}{10.8}{\rmdefault}{\mddefault}{\updefault}$1|0$}}}
\put(1082,980){\makebox(0,0)[lb]{{\SetFigFont{9}{10.8}{\rmdefault}{\mddefault}{\updefault}$1|0$}}}
\put(2584,2925){\makebox(0,0)[lb]{{\SetFigFont{9}{10.8}{\rmdefault}{\mddefault}{\updefault}$0|1$}}}
\put(2584,2700){\makebox(0,0)[lb]{{\SetFigFont{9}{10.8}{\rmdefault}{\mddefault}{\updefault}$1|0$}}}
\put(1609,1500){\makebox(0,0)[lb]{{\SetFigFont{9}{10.8}{\rmdefault}{\mddefault}{\updefault}$0|1$}}}
\put(3409,75){\makebox(0,0)[lb]{{\SetFigFont{9}{10.8}{\rmdefault}{\mddefault}{\updefault}$d_1$}}}
\put(4609,1500){\makebox(0,0)[lb]{{\SetFigFont{9}{10.8}{\rmdefault}{\mddefault}{\updefault}$d_2$}}}
\put(3409,2925){\makebox(0,0)[lb]{{\SetFigFont{9}{10.8}{\rmdefault}{\mddefault}{\updefault}$d_3$}}}
\end{picture}
\end{center}
\caption{A six-state automaton from $\mathfrak F$\label{sixstate}}
\end{figure}

If $\Ac= (Q,\{0,1\},\delta,\lambda)$ is an invertible automaton over
a binary alphabet, we say that a state $q\in Q$ is \emph{active} if
$q(0)=1$, $q(1)=0$; otherwise we say that $q$ is \emph{inactive}. We
consider a family $\mathfrak F$ of bireversible automata over a
binary alphabet defined as follows.  The state set of a member
$\Ac\in \mathfrak F$ is $Q_n = \{a,b,c,d_1,d_2,\ldots,d_n\}$, where
$n\geq 1$ can be any odd number. In wreath product coordinates
(c.f.~\eqref{wreathcoords}), $\Ac$ is given by:
\begin{gather*}
a=(01)(c,b),\ b=(01)(b,c),\ c=\sigma_0(d_1,d_1)\\ d_i =
\sigma_i(d_{i+1},d_{i+1}), 1\leq i\leq n-1\ \text{and}\ d_n =
\sigma_n(a,a)
\end{gather*}
 where the only restrictions on the $\sigma_i\in
S_{\{0,1\}}$, $i=0,\ldots,n$, is that an {\em odd} number of them
are not the identity
--- that is, an odd number of the states $c,d_1,\ldots,d_n$ are active.
In Figures~\ref{fourstate} and~\ref{sixstate}, respectively, we give
four-state and
six-state examples from the family $\mathfrak F$.

The inverse of $\Ac$ is given in wreath product coordinates
\eqref{wreathcoords} by:
\begin{gather*}
a\inv=(01)(b\inv,c\inv),\ b\inv=(01)(c\inv,b\inv),\
c\inv=\sigma_0(d_1\inv,d_1\inv)\\
d_i\inv = \sigma_i(d_{i+1}\inv,d_{i+1}\inv), 1\leq i\leq n-1\
\text{and}\ d_n\inv = \sigma_n(a\inv,a\inv)
\end{gather*}

\begin{proposition}\label{isbireversible}
Any automaton from the family $\mathfrak{F}$ is bireversible.
\end{proposition}
\begin{proof}
Let $\Ac$ have state set $Q_n$.   Then the input letter $0$ acts on
the states of $\Ac$ as the $(n+2)$-cycle $(acd_1\cdots d_n)$ while the
input letter $1$ acts on the states of $\Ac$ as the $(n+3)$-cycle
$(abcd_1\cdots d_n)$.  Thus $\Ac$ is reversible.  Similarly, $0$
acts on the states of $\Ac\inv$ as the $(n+3)$-cycle $(a\inv b\inv
c\inv d_1\inv \cdots d_n\inv)$ and $1$ acts as the $(n+2)$-cycle
$(a\inv c\inv d_1\inv \cdots d_n\inv)$, so $\Ac\inv$ is reversible
and hence $\Ac$ is bireversible.
\end{proof}

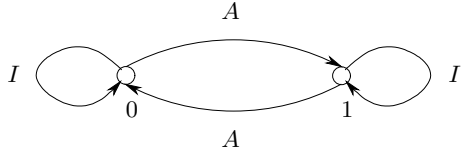
\begin{figure}[htbp]
\begin{center}
\setlength{\unitlength}{0.00062500in}
\begingroup\makeatletter\ifx\SetFigFont\undefined%
\gdef\SetFigFont#1#2#3#4#5{%
  \reset@font\fontsize{#1}{#2pt}%
  \fontfamily{#3}\fontseries{#4}\fontshape{#5}%
  \selectfont}%
\fi\endgroup%
\begin{picture}(4892,1575)(0,-10)
\put(2400.000,-787.500){\arc{3825.000}{4.2224}{5.2023}}
\blacken\path(3179.145,926.343)(3300.000,900.000)(3205.760,980.117)(3224.717,937.261)(3179.145,926.343)
\put(2400.000,2437.000){\arc{3824.118}{1.0807}{2.0609}}
\blacken\path(1620.853,723.645)(1500.000,750.000)(1594.233,669.874)(1575.280,712.732)(1620.853,723.645)
\put(1500,825){\ellipse{150}{150}}
\put(3303,823){\ellipse{150}{150}}
\path(1464,872)(1463,873)(1459,877)
    (1454,882)(1445,890)(1433,901)
    (1419,914)(1402,929)(1383,944)
    (1363,961)(1341,977)(1318,993)
    (1294,1008)(1269,1023)(1242,1035)
    (1213,1046)(1183,1055)(1150,1062)
    (1116,1065)(1080,1065)(1048,1061)
    (1018,1055)(990,1047)(965,1039)
    (943,1030)(923,1022)(906,1013)
    (891,1005)(878,998)(866,990)
    (854,982)(843,975)(832,966)
    (822,957)(810,947)(799,935)
    (788,921)(777,906)(766,888)
    (758,869)(752,847)(750,825)
    (753,803)(760,781)(770,761)
    (782,743)(794,726)(806,712)
    (818,699)(830,688)(842,678)
    (853,668)(865,659)(877,651)
    (890,642)(904,634)(919,625)
    (937,615)(957,606)(979,596)
    (1005,587)(1033,579)(1063,572)
    (1095,568)(1130,568)(1164,571)
    (1195,578)(1225,588)(1252,600)
    (1277,613)(1301,628)(1323,644)
    (1344,662)(1364,679)(1383,697)
    (1400,714)(1415,729)(1428,743)
    (1439,755)(1456,774)
\blacken\path(1398.342,664.567)(1456.000,774.000)(1353.627,704.575)(1399.989,711.400)(1398.342,664.567)
\path(3336,872)(3337,873)(3341,877)
    (3346,882)(3355,890)(3367,901)
    (3381,914)(3398,929)(3417,944)
    (3437,961)(3459,977)(3482,993)
    (3506,1008)(3531,1023)(3558,1035)
    (3587,1046)(3617,1055)(3650,1062)
    (3684,1065)(3720,1065)(3752,1061)
    (3782,1055)(3810,1047)(3835,1039)
    (3857,1030)(3877,1022)(3894,1013)
    (3909,1005)(3922,998)(3934,990)
    (3946,982)(3957,975)(3968,966)
    (3978,957)(3990,947)(4001,935)
    (4012,921)(4023,906)(4034,888)
    (4042,869)(4048,847)(4050,825)
    (4047,803)(4040,781)(4030,761)
    (4018,743)(4006,726)(3994,712)
    (3982,699)(3970,688)(3958,678)
    (3947,668)(3935,659)(3923,651)
    (3910,642)(3896,634)(3881,625)
    (3863,615)(3843,606)(3821,596)
    (3795,587)(3767,579)(3737,572)
    (3705,568)(3670,568)(3636,571)
    (3605,578)(3575,588)(3548,600)
    (3523,613)(3499,628)(3477,644)
    (3456,662)(3436,679)(3417,697)
    (3400,714)(3385,729)(3372,743)
    (3361,755)(3344,774)
\blacken\path(3446.373,704.575)(3344.000,774.000)(3401.658,664.567)(3400.011,711.400)(3446.373,704.575)
\put(3300,475){\makebox(0,0)[lb]{{\SetFigFont{9}{10.8}{\rmdefault}{\mddefault}{\updefault}$1$}}}
\put(1500,475){\makebox(0,0)[lb]{{\SetFigFont{9}{10.8}{\rmdefault}{\mddefault}{\updefault}$0$}}}
\put(2300,225){\makebox(0,0)[lb]{{\SetFigFont{9}{10.8}{\rmdefault}{\mddefault}{\updefault}$A$}}}
\put(4200,775){\makebox(0,0)[lb]{{\SetFigFont{9}{10.8}{\rmdefault}{\mddefault}{\updefault}$I$}}}
\put(2300,1300){\makebox(0,0)[lb]{{\SetFigFont{9}{10.8}{\rmdefault}{\mddefault}{\updefault}$A$}}}
\put(510,775){\makebox(0,0)[lb]{{\SetFigFont{9}{10.8}{\rmdefault}{\mddefault}{\updefault}$I$}}}
\end{picture}
\end{center}
\caption{The transition diagram for the right-to-left scanning automata $\Dc$,
$\Ec$ and $\Fc$\label{automatonD}}
\end{figure}

We fix for the rest of this section a member
$\Ac=(Q_n,\{0,1\},\delta,\lambda)$ of the family $\mathfrak{F}$. Let
$A$ denote the set of active states of $\Ac^{\pm}$ and $I$ the set
of inactive states.  Notice that a state $q$ is active if and only
if $q\inv$ is active; that is $A$ and $I$ are closed under the
involution.   Let $\Dc = \wh {\Ac^{\pm}}$ be the dual of the
disjoint union of $\Ac$ and $\Ac\inv$.  The state set of $\Dc$ is
$\{0,1\}$. Since active states of $\Ac$ switch $0$ and $1$, while
inactive states do not, it follows that each input letter from $A$ to
$\Dc$ switches states, while each input letter from $I$ does not.
Thus the transition diagram of $\Dc$ has the form in
Figure~\ref{automatonD}.

The proof of Proposition~\ref{isbireversible} shows that the output
function of $\Dc_0$ on letters is given by the permutation
\[\delta_0=(acd_1\cdots d_n)(a\inv b\inv c\inv d_1\inv\cdots d_n\inv)\] while
the output function of $\Dc_1$ on letters is given by the
permutation \[\delta_1=(abcd_1\cdots d_n)(a\inv c\inv d_1\inv\cdots
d_n\inv)\]

In wreath product coordinates \eqref{wreathcoordsdual}, and using
functional notation for $Q_n^{\pm}$-tuples, we can then write $\Dc_0
= (D_0,\delta_0)$, $\Dc_1 = (D_1,\delta_1)$ where
$D_0,D_1:Q_n^\pm\to \GG(\Dc)$ are given by
\[qD_i =
\begin{cases} \Dc_{\ov i} & q\in A\\ \Dc_i & q\in I\end{cases}\]
with the notation $\ov 0=1$, $\ov 1=0$; this notation shall reoccur
throughout the rest of the paper without comment.

The right-to-left scanning automaton $\Dc$ has the following useful property,
which the reader easily checks: if
one reverses all the arrows in the Moore diagram of $\Dc$, then you
get back exactly $\Dc$, only with the states $0$ and $1$ reversed.
Basically, the arrows with left side labelled by elements of $I$ remain
as they were,
while the arrows with left side labelled by $A$ from $0$ and $1$ are
switched with the corresponding arrows from $1$ to $0$.
  One can then obtain via an easy induction that if $u$ is in
the orbit under $\Sc(\Dc)$ of $v$, then $u^{\rho}$ is in the orbit
under $\Sc(\Dc)$ of $v^{\rho}$.  But for invertible automata, the
orbit of the group and the semigroup are the same by
Lemma~\ref{semigrouporbits}, so we have proven:

\begin{lemma}\label{canreverse}
Let $u,v\in Q_n^{\pm}$.  Then $u$ and $v$ are in the same
$\GG(\Dc)$-orbit if and only if $u^{\rho}$ and $v^{\rho}$ are in the
same $\GG(\Dc)$-orbit.
\end{lemma}

\subsection{The supporting cast}
We shall also need to define and study several auxiliary
right-to-left scanning automata that will play a key role in proving
transitivity of $\GG(\Dc)$ on patterns. Since we shall deal
exclusively for the next few sections with right-to-left scanning
automata, we shall use the convention that permutations act on the
right of their arguments and compose them appropriately: e.g.\
$(ab)(abc) = (ac)$. Define a right-to-left scanning automaton $\Ec =
(\{0,1\},Q_n^{\pm},\delta_E,\lambda_E)$ using the same transition
diagram as $\Dc$ (see Figure~\ref{automatonD}) but with a different
output function, modelled on $\mathsf E$ from Figure~\ref{aleshinE}.
Let $\varepsilon_0 = (a\inv b\inv)$ and $\varepsilon_1 = (ab)$. Then
we have $\Ec_0$ act on letters by $\varepsilon_0$ and $\Ec_1$ by
$\varepsilon_1$.  So in wreath product coordinates
\eqref{wreathcoordsdual}, $\Ec_0 = (E_0,\varepsilon_0)$,
$\Ec_1=(E_1,\varepsilon_1)$ where \[qE_i =
\begin{cases} \Ec_{\ov i} & q\in A\\ \Ec_i & q\in I\end{cases}\]

 Recall that if $\sigma\in
S_{Q_n}$, then $\ov\sigma$ denotes the induced involution-preserving
permutation of $S_{Q_n^{\pm}}$. One can verify by direct
computation:
\begin{align*}
\varepsilon_0 \ov {(acd_1\cdots d_n)} &= (a\inv b\inv)\ov
{(acd_1\cdots
d_n)} = \delta_0\\
\varepsilon_1 \ov {(acd_1\cdots d_n)} &= (ab)\ov {(acd_1\cdots
d_n)} = \delta_1\\
\varepsilon_0 \ov {(abcd_1\cdots d_n)} &= (a\inv b\inv)\ov
{(abcd_1\cdots
d_n)} = \delta_1\\
\varepsilon_1 \ov {(abcd_1\cdots d_n)} &= (ab)\ov {(abcd_1\cdots
d_n)} = \delta_0
\end{align*}
We can then deduce from the right-to-left scanning dual to
Proposition~\ref{composewithauto} the following lemma.

\begin{lemma}\label{DvsE}
We have the equalities:
\begin{enumerate}
\item  $\Dc_0=\Ec_0(acd_1\cdots d_n)^{\pm}$\item  $\Dc_1 =\Ec_1(acd_1\cdots
d_n)^{\pm}$\item  $\Dc_0 = \Ec_1(abcd_1\cdots d_n)^{\pm}$ \item
$\Dc_1 =\Ec_0(abcd_1\cdots d_n)^{\pm}$
\end{enumerate}
\end{lemma}

Lemma~\ref{DvsE} indicates that we should study $\GG(\Ec)$.  The
analysis is completely analogous to the situation for $\mathsf E$
from~\cite{free1}.

\begin{lemma}\label{studyE}
The group $\GG (\Ec)$ is a Klein $4$-group.  More precisely:
\begin{enumerate}
\item $\Ec_0^2=\Ec_1^2=1$
\item $\Ec_i(ab)^\pm = \Ec_{\ov i}$, $i=0,1$
\item $\Ec_1\Ec_0 = (ab)^{\pm} = \Ec_0\Ec_1$
\end{enumerate}
\end{lemma}
\begin{proof}
First observe that $\Ec = \Ec\inv$.  Indeed the only edges of the
Moore diagram for which both sides of the label are not the same
are: $0\arr{a\inv \mid b\inv}1$, $0\arr {b\inv\mid a\inv}1$,
$1\arr{a\mid b}0$ and $1\arr{b\mid a} 0$.  But interchanging the
left and right sides of the labels just switches the first two
arrows and switches the last two arrows thereby having no effect on
the automaton.  It follows $\Ec_0=\Ec_0\inv$ and $\Ec_1=\Ec_1\inv$,
establishing (1).

For (2) and (3), observe that $(a\inv b\inv)\ov{(ab)} = (ab)$ and
$(ab)\ov{(ab)} = (a\inv b\inv)$.  Thus
Proposition~\ref{composewithauto} shows that $\Ec_0 =
\Ec_1(ab)^{\pm}$ and $\Ec_1= \Ec_0(ab)^{\pm}$.  Thus \[\Ec_0\Ec_1 =
\Ec_0^2(ab)^{\pm} = (ab)^{\pm}\]  Hence \[ \Ec_1\Ec_0 =
(\Ec_0\Ec_1)\inv = ((ab)^{\pm})\inv = (ab)^{\pm}\] finishing the
proof.
\end{proof}

As a consequence, we obtain:
\begin{lemma}\label{getabbc}
One has:
\begin{enumerate}
\item $\Dc_0\inv \Dc_1 = (bc)^{\pm}$
\item $\Dc_0\Dc_1\inv = (ab)^{\pm}$
\end{enumerate}
Hence, for all $\sigma\in S_{a,b,c}$, one has $\sigma^{\pm}\in
\GG(\Dc)$.
\end{lemma}
\begin{proof}
We calculate:
\begin{align*}
\Dc_0\inv \Dc_1 &= (d_n\cdots d_1ca)^{\pm}\Ec_0\Ec_0(abcd_1\cdots
d_n)^{\pm} = (bc)^{\pm}\\
\Dc_0\Dc_1\inv &= \Ec_0(acd_1\cdots d_n)^{\pm}(d_n\cdots
d_1cba)^{\pm}\Ec_0 =  \Ec_0(ab)^{\pm}\Ec_0 = (ab)^{\pm}
\end{align*}
where the last equality follows from Lemma~\ref{studyE}, which
implies that $\Ec_0$ and $(ab)^{\pm}$ commute.
\end{proof}

Our next auxiliary right-to-left scanning automaton, called $\Fc$, again has
transition diagram from Figure~\ref{automatonD}.  Let $\varphi_0 =
(a\inv b\inv)\ov{(cd_1\cdots d_n)}$ and $\varphi_1 =
(ab)\ov{(cd_1\cdots d_n)}$.  Then the output of $\Fc_0$ on letters
is given by the permutation $\varphi_0$ and the output of $\Fc_1$ on
letters is given by the permutation $\varphi_1$.  So in wreath
product coordinates \eqref{wreathcoordsdual},
$\Fc_0=(F_0,\varphi_0)$ and $\Fc_1 = (F_1,\varphi_1)$ where \[qF_i =
\begin{cases} \Fc_{\ov i} & q\in A\\ \Fc_i & q\in I\end{cases}\]
Proposition~\ref{composewithauto} shows that
\[\Fc_0=\Ec_0(cd_1\cdots d_n)^\pm = \Dc_0(ac)^\pm\ \text{and}\
\Fc_1=\Ec_1(cd_1\cdots d_n)^\pm = \Dc_1(ac)^\pm\] where the
equalities involving $\Dc$ come from Lemma~\ref{DvsE}.  As a
consequence of Lemma~\ref{getabbc}, we obtain:

\begin{lemma}\label{Fin}
$\Fc_0,\Fc_1\in \GG(\Dc)$.  In fact, $\GG(\Dc) = \langle
\Fc_0,\Fc_1, (ac)^{\pm}\rangle$.
\end{lemma}

This leads us to investigate further the structure of $\Fc$.

\begin{lemma}\label{studyF}
The following facts hold for $\Fc$:
\begin{enumerate}
\item $\Fc_i = \Fc_{\ov i}(ab)^\pm = (ab)^\pm\Fc_{\ov i}$, $i=0,1$
\item $\Fc_0\Fc_1=\Fc_1\Fc_0$
\item if $w(x,y)\in \{x,y\}^*$, then \[w(\Fc_0,\Fc_1) =\begin{cases} \Fc_0^{|w|} & \text{if $w$ has
an even number of $y$s}\\ \Fc_0^{|w|}(ab)^\pm = \Fc_0^{|w|-1}\Fc_1 &
\text{if $w$ has an odd number of $y$s}\end{cases}\]
\end{enumerate}
\end{lemma}
\begin{proof}
Let $i\in \{0,1\}$.  Now, Lemma~\ref{studyE} gives us:
\[\Fc_i = \Ec_i(cd_1\cdots d_n)^\pm = \Ec_{\ov i}(ab)^\pm(cd_1\cdots
d_n)^\pm = \Fc_{\ov i}(ab)^\pm\] But Lemma~\ref{studyE} implies that
$(ab)^\pm$ commutes with $\Ec_{\ov i}$, so we also get: \[\Fc_i
=\Ec_{\ov i}(ab)^\pm(cd_1\cdots d_n)^\pm=(ab)^\pm\Ec_{\ov
i}(cd_1\cdots d_n)^\pm =(ab)^\pm\Fc_{\ov i}\] For (2), we now have
$\Fc_0\Fc_1 = \Fc_0(ab)^{\pm}\Fc_0 = \Fc_1\Fc_0$. Item (3) follows
immediately (1) and (2).
\end{proof}

Now if $m\geq 1$, then the sections of $\Fc_0^m$ are of the form
$w(\Fc_0,\Fc_1)$ where $|w|=m$.  Thus $\Fc_0^m$ has at most two
sections, itself $\Fc_0^m$ and $\Fc_0^m(ab)^\pm$.

\begin{figure}[htbp]
\begin{center}
\setlength{\unitlength}{0.00062500in}
\begingroup\makeatletter\ifx\SetFigFont\undefined%
\gdef\SetFigFont#1#2#3#4#5{%
  \reset@font\fontsize{#1}{#2pt}%
  \fontfamily{#3}\fontseries{#4}\fontshape{#5}%
  \selectfont}%
\fi\endgroup%
\begin{picture}(4892,1575)(0,-10)
\put(2400.000,-787.500){\arc{3825.000}{4.2224}{5.2023}}
\blacken\path(3179.145,926.343)(3300.000,900.000)(3205.760,980.117)(3224.717,937.261)(3179.145,926.343)
\put(2400.000,2437.000){\arc{3824.118}{1.0807}{2.0609}}
\blacken\path(1620.853,723.645)(1500.000,750.000)(1594.233,669.874)(1575.280,712.732)(1620.853,723.645)
\put(1500,825){\ellipse{150}{150}}
\put(3303,823){\ellipse{150}{150}}
\path(1464,872)(1463,873)(1459,877)
    (1454,882)(1445,890)(1433,901)
    (1419,914)(1402,929)(1383,944)
    (1363,961)(1341,977)(1318,993)
    (1294,1008)(1269,1023)(1242,1035)
    (1213,1046)(1183,1055)(1150,1062)
    (1116,1065)(1080,1065)(1048,1061)
    (1018,1055)(990,1047)(965,1039)
    (943,1030)(923,1022)(906,1013)
    (891,1005)(878,998)(866,990)
    (854,982)(843,975)(832,966)
    (822,957)(810,947)(799,935)
    (788,921)(777,906)(766,888)
    (758,869)(752,847)(750,825)
    (753,803)(760,781)(770,761)
    (782,743)(794,726)(806,712)
    (818,699)(830,688)(842,678)
    (853,668)(865,659)(877,651)
    (890,642)(904,634)(919,625)
    (937,615)(957,606)(979,596)
    (1005,587)(1033,579)(1063,572)
    (1095,568)(1130,568)(1164,571)
    (1195,578)(1225,588)(1252,600)
    (1277,613)(1301,628)(1323,644)
    (1344,662)(1364,679)(1383,697)
    (1400,714)(1415,729)(1428,743)
    (1439,755)(1456,774)
\blacken\path(1398.342,664.567)(1456.000,774.000)(1353.627,704.575)(1399.989,711.400)(1398.342,664.567)
\path(3336,872)(3337,873)(3341,877)
    (3346,882)(3355,890)(3367,901)
    (3381,914)(3398,929)(3417,944)
    (3437,961)(3459,977)(3482,993)
    (3506,1008)(3531,1023)(3558,1035)
    (3587,1046)(3617,1055)(3650,1062)
    (3684,1065)(3720,1065)(3752,1061)
    (3782,1055)(3810,1047)(3835,1039)
    (3857,1030)(3877,1022)(3894,1013)
    (3909,1005)(3922,998)(3934,990)
    (3946,982)(3957,975)(3968,966)
    (3978,957)(3990,947)(4001,935)
    (4012,921)(4023,906)(4034,888)
    (4042,869)(4048,847)(4050,825)
    (4047,803)(4040,781)(4030,761)
    (4018,743)(4006,726)(3994,712)
    (3982,699)(3970,688)(3958,678)
    (3947,668)(3935,659)(3923,651)
    (3910,642)(3896,634)(3881,625)
    (3863,615)(3843,606)(3821,596)
    (3795,587)(3767,579)(3737,572)
    (3705,568)(3670,568)(3636,571)
    (3605,578)(3575,588)(3548,600)
    (3523,613)(3499,628)(3477,644)
    (3456,662)(3436,679)(3417,697)
    (3400,714)(3385,729)(3372,743)
    (3361,755)(3344,774)
\blacken\path(3446.373,704.575)(3344.000,774.000)(3401.658,664.567)(3400.011,711.400)(3446.373,704.575)
\put(3300,475){\makebox(0,0)[lb]{{\SetFigFont{9}{10.8}{\rmdefault}{\mddefault}{\updefault}$1$}}}
\put(1500,475){\makebox(0,0)[lb]{{\SetFigFont{9}{10.8}{\rmdefault}{\mddefault}{\updefault}$0$}}}
\put(4200,900){\makebox(0,0)[lb]{{\SetFigFont{9}{10.8}{\rmdefault}{\mddefault}{\updefault}$a|b\quad
a\inv|b\inv$}}}
\put(4200,675){\makebox(0,0)[lb]{{\SetFigFont{9}{10.8}{\rmdefault}{\mddefault}{\updefault}$b|a\quad
b\inv|a\inv$}}}
\put(2000,1300){\makebox(0,0)[lb]{{\SetFigFont{9}{10.8}{\rmdefault}{\mddefault}{\updefault}$Q^\pm\setminus
\{a,b\}^\pm$}}}
\put(2000,225){\makebox(0,0)[lb]{{\SetFigFont{9}{10.8}{\rmdefault}{\mddefault}{\updefault}$Q^\pm\setminus
\{a,b\}^\pm$}}}
\put(-500,675){\makebox(0,0)[lb]{{\SetFigFont{9}{10.8}{\rmdefault}{\mddefault}{\updefault}$b|b\quad
b\inv|b\inv$}}}
\put(-500,900){\makebox(0,0)[lb]{{\SetFigFont{9}{10.8}{\rmdefault}{\mddefault}{\updefault}$a|a\quad
a\inv|a\inv$}}}
\end{picture}
\end{center}
\caption{The Moore diagram for the right-to-left scanning automaton
$\Gc$\label{automatonG}}
\end{figure}
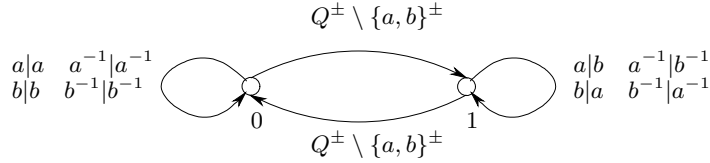

Let $n\geq 1$ be the odd number associated to our automaton $\Ac$
with state set $Q_n$ from the family $\mathfrak{F}$.  We want to
show that $\Fc_0^{n+1}$ has two sections.  More precisely let
$\mathcal G$ be the right-to-left scanning automaton with Moore diagram in
Figure~\ref{automatonG}, where we take the convention that if the
right side of the label of an edge equals the left side then we
omit the right side.

In wreath product coordinates \eqref{wreathcoordsdual}, we have
$\Gc_0 = (G_0,1)$, $\Gc_1 = (G_1,\ov{(ab)})$ where \[qG_i =
\begin{cases} G_i & q\in \{a,b\}^{\pm}\\ G_{\ov i} &
\text{else}\end{cases}\]  The proof of our next lemma is essentially
an exercise is computing powers of elements in wreath products.  It
is the key place where we use that $\Ac$ has an odd number of active
states.

\begin{lemma}\label{StudyG}
$\Fc^{n+1}_0 = \Gc_0$, $\Fc_0^{n+1}(ab)^\pm = \Fc_0^n\Fc_1 = \Gc_1$.
In particular, $\Gc_0,\Gc_1\in \GG(\Dc)$.
\end{lemma}
\begin{proof}
Recall that in wreath product coordinates, we have $\Fc_0 =
(F_0,\varphi_0)$ where $F_0$ takes active states to $\Fc_1$ and
inactive states to $\Fc_0$ and $\varphi_0 =(a\inv
b\inv)\ov{(cd_1\cdots d_n)}$.  Then $\varphi_0^{n+1} =1$, as $n+1$
is even.  So $\Fc_0^{n+1} = (f,1)$ where $f =
F_0{}^{\varphi_0}\!{F_0}{}^{\varphi_0^2}\!{F_0}\cdots
{}^{\varphi_0^{n}}\!{F_0}$. Notice that $\{a,b\}^{\pm}$ is invariant
under $\varphi_0$ and that $F_0$ takes on the constant value $\Fc_1$
on $\{a,b\}^{\pm}$, since these states are all active. It follows
that $qf = \Fc_1^{n+1} = \Fc_0^{n+1}$ (this last equality by
Lemma~\ref{studyF}) for $q\in \{a,b\}^{\pm}$. On the other hand,
$\{c,d_1,\ldots, d_n\}$ and $\{c\inv, d_1\inv,\ldots, d_n\inv\}$ are
cycles of length $n+1$ for $\varphi_0$. Thus, for $q\in Q\setminus
\{a,b\}$ and $e=\pm 1$, one has that $q^ef = \prod_{r\in Q} r^eF_0$
where the product is taken in cyclic order starting from $q$.
Actually, the order doesn't matter since $F_0$ only takes on the
values $\Fc_0$ and $\Fc_1$ and these elements commute. So one has
$q^ef=w(\Fc_0,\Fc_1)$ where $w(x,y)$ is a word of length $n+1$. The
occurrences of the variable $x$ come from inactive states in
$Q\setminus \{a,b\}$ while the occurrences of the variable $y$ come
from the active states in this set.  Since by hypothesis $\Ac$ has
an odd number of active states in this set, we conclude by
Lemma~\ref{studyF} that $q^ef = \Fc_0^{n+1}(ab)^{\pm} =
\Fc_0^n\Fc_1$.

Let us abuse notation and, for $g\in \GG(\Dc)$, denote also by $g$
the function $\GG(\Dc)^{Q^\pm}$ that takes on the constant value
$g$. Then one has that \[\Fc_0^n\Fc_1=\Fc_0^{n+1}(ab)^{\pm} =
(f,1)((ab)^{\pm},\ov{(ab)}) = (f(ab)^{\pm},\ov{(ab)})\]  In
particular, Lemma~\ref{studyF} implies that if $qf=\Fc_i$, then
$qf(ab)^\pm = \Fc_{\overline i}$. Converting wreath product
coordinates to automata, we see that on inputs in $\{a,b\}^{\pm}$,
both $\Fc_0^{n+1}$ and $\Fc_0^n\Fc_1$ remain in the same state,
while on inputs from $Q^\pm\setminus\{a,b\}^{\pm}$ they both switch
states.  On letters, $\Fc_0^{n+1}$ acts as the identity, while
$\Fc_0^n\Fc_1$ acts as the permutation $\ov{(ab)}$. Thus $\Fc_0
=\Gc_0$, $\Fc_1=\Gc_1$, as required.
\end{proof}

We remark that $\Gc=\Gc\inv$, and so the order of $\Fc_0$ (and hence
$\Fc_1$) is $2(n+1)$.  However, we shall not use this fact anywhere.

\subsection{A second family of automata}
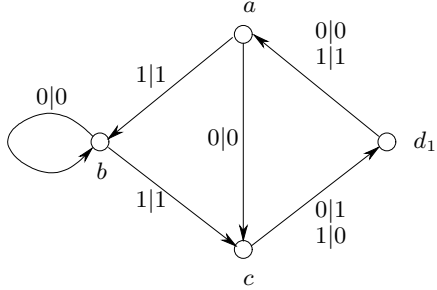
\begin{figure}[htbp]
\begin{center}
\setlength{\unitlength}{0.00062500in}
\begingroup\makeatletter\ifx\SetFigFont\undefined%
\gdef\SetFigFont#1#2#3#4#5{%
  \reset@font\fontsize{#1}{#2pt}%
  \fontfamily{#3}\fontseries{#4}\fontshape{#5}%
  \selectfont}%
\fi\endgroup%
\begin{picture}(3593,2395)(0,-10)
\put(784,1200){\ellipse{150}{150}}
\put(1984,2100){\ellipse{150}{150}}
\put(1984,300){\ellipse{150}{150}}
\put(3184,1200){\ellipse{150}{150}} \path(856,1152)(1921,327)
\blacken\path(1807.762,376.771)(1921.000,327.000)(1844.506,424.204)(1854.594,378.441)(1807.762,376.771)
\path(1893,2075)(853,1227)
\blacken\path(927.044,1326.083)(853.000,1227.000)(964.960,1279.582)(918.102,1280.083)(927.044,1326.083)
\path(2051,330)(3116,1155)
\blacken\path(3039.506,1057.796)(3116.000,1155.000)(3002.762,1105.229)(3049.594,1103.559)(3039.506,1057.796)
\path(3115,1247)(2075,2095)
\blacken\path(2186.960,2042.418)(2075.000,2095.000)(2149.044,1995.917)(2140.102,2041.917)(2186.960,2042.418)
\path(1984,2025)(1984,375)
\blacken\path(1954.000,495.000)(1984.000,375.000)(2014.000,495.000)(1984.000,459.000)(1954.000,495.000)
\path(726,1246)(725,1247)(721,1251)
    (716,1256)(707,1264)(695,1275)
    (681,1288)(664,1303)(645,1318)
    (625,1335)(603,1351)(580,1367)
    (556,1382)(531,1397)(504,1409)
    (475,1420)(445,1429)(412,1436)
    (378,1439)(342,1439)(310,1435)
    (280,1429)(252,1421)(227,1413)
    (205,1404)(185,1396)(168,1387)
    (153,1379)(140,1372)(128,1364)
    (116,1356)(105,1349)(94,1340)
    (84,1331)(72,1321)(61,1309)
    (50,1295)(39,1280)(28,1262)
    (20,1243)(14,1221)(12,1199)
    (15,1177)(22,1155)(32,1135)
    (44,1117)(56,1100)(68,1086)
    (80,1073)(92,1062)(104,1052)
    (115,1042)(127,1033)(139,1025)
    (152,1016)(166,1008)(181,999)
    (199,989)(219,980)(241,970)
    (267,961)(295,953)(325,946)
    (357,942)(392,942)(426,945)
    (457,952)(487,962)(514,974)
    (539,987)(563,1002)(585,1018)
    (606,1036)(626,1053)(645,1071)
    (662,1088)(677,1103)(690,1117)
    (701,1129)(718,1148)
\blacken\path(660.342,1038.567)(718.000,1148.000)(615.627,1078.575)(661.989,1085.400)(660.342,1038.567)
\put(1984,0){\makebox(0,0)[lb]{{\SetFigFont{9}{10.8}{\rmdefault}{\mddefault}{\updefault}$c$}}}
\put(754,890){\makebox(0,0)[lb]{{\SetFigFont{9}{10.8}{\rmdefault}{\mddefault}{\updefault}$b$}}}
\put(1984,2290){\makebox(0,0)[lb]{{\SetFigFont{9}{10.8}{\rmdefault}{\mddefault}{\updefault}$a$}}}
\put(252,1467){\makebox(0,0)[lb]{{\SetFigFont{9}{10.8}{\rmdefault}{\mddefault}{\updefault}$0|0$}}}
\put(1082,1655){\makebox(0,0)[lb]{{\SetFigFont{9}{10.8}{\rmdefault}{\mddefault}{\updefault}$1|1$}}}
\put(1082,605){\makebox(0,0)[lb]{{\SetFigFont{9}{10.8}{\rmdefault}{\mddefault}{\updefault}$1|1$}}}
\put(2584,525){\makebox(0,0)[lb]{{\SetFigFont{9}{10.8}{\rmdefault}{\mddefault}{\updefault}$0|1$}}}
\put(2584,300){\makebox(0,0)[lb]{{\SetFigFont{9}{10.8}{\rmdefault}{\mddefault}{\updefault}$1|0$}}}
\put(2584,2025){\makebox(0,0)[lb]{{\SetFigFont{9}{10.8}{\rmdefault}{\mddefault}{\updefault}$0|0$}}}
\put(2584,1800){\makebox(0,0)[lb]{{\SetFigFont{9}{10.8}{\rmdefault}{\mddefault}{\updefault}$1|1$}}}
\put(1684,1125){\makebox(0,0)[lb]{{\SetFigFont{9}{10.8}{\rmdefault}{\mddefault}{\updefault}$0|0$}}}
\put(3409,1125){\makebox(0,0)[lb]{{\SetFigFont{9}{10.8}{\rmdefault}{\mddefault}{\updefault}$d_1$}}}
\end{picture}
\end{center}
\caption{A four-state automaton from $\mathfrak
F'$\label{fourstateprime}}
\end{figure}

Let $\mathfrak F'$ be the family of all automata that can be
obtained from a member of $\mathfrak F$ by switching which states
are active  and  inactive.  So, for example, the Moore diagram of
the element of $\mathfrak F'$ corresponding to the automaton in
Figure~\ref{fourstate} is given by Figure~\ref{fourstateprime}.

If $\Ac=(Q_n,\{0,1\},\delta,\lambda)\in \mathfrak F$ and $\Ac'$ is
the corresponding automaton in $\mathfrak F'$, then in the
terminology of Proposition~\ref{composewithauto}, $\Ac' = (01)[\Ac]$
and so $\GG(\Ac') = \langle (01)^*\Ac_q\mid q\in Q_n\}$.  In wreath
product coordinates \eqref{wreathcoords}, $\Ac'$ is given by:
\begin{gather*} a=(c,b),\ b=(b,c),\ c=\sigma_0'(d_1,d_1)\\
d_i = \sigma_i'(d_{i+1},d_{i+1}), 1\leq i\leq n-1\ \text{and}\ d_n =
\sigma_n'(a,a)
\end{gather*}
where $\sigma_i'$ is the unique element of $S_2\setminus
\{\sigma_i\}$.  Thus $(\Ac')\inv$ is described in wreath product
coordinates by:
\begin{gather*} a\inv=(c\inv,b\inv),\ b\inv=(b\inv,c\inv),\ c\inv=\sigma_0'(d_1\inv,d_1\inv)\\
d_i\inv = \sigma_i'(d_{i+1}\inv,d_{i+1}\inv), 1\leq i\leq n-1\
\text{and}\ d_n\inv = \sigma_n'(a\inv,a\inv)
\end{gather*}
from which it is immediate that $(\Ac')\inv =\Ac'$ and so each
element $\Ac'_q$, with $q\in Q_n$, is its own inverse.  We summarize
this discussion as a lemma.

\begin{lemma}\label{secondfamilydetails}
Let $\Ac=(Q,A,\delta,\lambda)\in \mathfrak F$ and let $(01)[\Ac]\in
\mathfrak F'$ be the corresponding automaton.  Then $(01)[\Ac]_q =
(01)^*\Ac_q$ and $((01)^*\Ac_q)^2=1$.
\end{lemma}

\section{Freeness results}
In this section, we prove, for $\Ac\in \mathfrak F$, that $\GG(\Ac)$
is a free group, freely generated by the states, and that
$\GG((01)[\Ac])$ is a free product of cyclic groups of order two,
again freely generated by the states.

\subsection{Freeness for automata in $\mathfrak F$}
 Let's fix $\Ac = (Q_n,\{0,1\},\delta,\lambda)$
from the family $\mathfrak F$ for the remainder of the section. Let
$\Dc,\Ec,\Fc,\Gc$ be as in the previous section.

\begin{lemma}\label{actpreservespatt}
The action of $\GG(\Dc)$ on $Q_n^\pm$ preserves patterns and sends
freely irreducible words to freely irreducible words.
\end{lemma}
\begin{proof}
Clearly $\GG(\Dc)$ preserves patterns since $\delta_0$ and
$\delta_1$ preserve $Q_n$ and $Q_n\inv$. To see that being freely
irreducible is preserved, it suffices to show that words with a
factor of the form $xx\inv$ with $x\in Q_n^\pm$ are preserved by the
action of $\GG(\Dc)$. Lemma~\ref{Fin} shows that $\GG(\Dc)$ is
generated by $\Fc_0$, $\Fc_1$ and $(ac)^\pm$. Clearly $(ac)^\pm$
preserves factors of the form $xx\inv$ $x\in Q_n^{\pm}$. On elements
of $((Q_n\setminus\{a,b\})^{\pm})^*$, both $\Fc_0$ and $\Fc_1$ acts
as $(cd_1\cdots d_m)^{\pm}$ and so preserve factors of the form
$xx\inv$ with $x\in(Q_n\setminus\{a,b\})^{\pm})$.  On the other
hand, $\Fc_0(xx\inv) = xx\inv$ for $x\in \{a,b\}$ while $\Fc_0(a\inv
a) = b\inv b$, $\Fc_0(b\inv b)=a\inv a$.  A similar computation
shows that $\Fc_1$ preserves such factors.  This completes the
proof.
\end{proof}

Next we show that each non-empty pattern contains an element that
changes the first letter of a word and hence a non-trivial element.

\begin{proposition}\label{actsnontrivlevel1}
Let $p\in \{\ast,\ast\inv\}^*$ be a non-empty pattern.   Then there
is a freely irreducible word $w\in (Q_n^{\pm})$, following the
pattern $p$, such that $w$ acts non-trivially on $\{0,1\}^*$.
\end{proposition}
\begin{proof}
We construct such a $w$ acting non-trivially on the first letter of
a word. Recall that an element $x^e$, with $x\in Q_n$, $e=\pm 1$,
acts non-trivially on the first letter of a word in $\{0,1\}^*$ if
and only if $x$ is an active state of $\Ac$.  First construct a word
$v$ by replacing each $\ast$ in $p$ by $a$ and each $\ast\inv$ by
$b\inv$. Evidently, $v$ is freely irreducible and follows $p$. If
$v$ acts non-trivially on the first letter of a word, we are done.
Else, since both $a$ and $b\inv$ are active, it follows that the
total number of $a$s and $b\inv$s is even. Choose any inactive state
$q$ from $\{c,d_1,\ldots,d_n\}$; since there are an even number of
elements in this set and by hypothesis an odd number are active, we
can find such a state $q$.  Now we obtain a new word $w$ by
replacing the first letter of $v$ by $q$ or $q\inv$ according to
whether the first letter of $p$ is $\ast$ or $\ast\inv$.  Then $w$
is still freely irreducible and now the number of letters in $w$
corresponding to actives states is odd.  So $w$ acts non-trivially
on the first letter of each word.
\end{proof}

We shall prove now that $\GG(\Dc)$ acts transitively on the set of
freely irreducible words following any given pattern.  Then
Corollary~\ref{freeness}, in light of
Proposition~\ref{actsnontrivlevel1}, will imply that $\GG(\Ac)$ is a
free group on $n+3$ generators, the elements $\Ac_q$, $q\in Q_n$,
being a free basis.

A key tool in proving the transitivity on patterns is a standard
lemma (c.f.~\cite{free1}) that plays an important role in the
induction argument. For $u\in (Q_n^\pm)^*$, we denote by $\St u$ the
stabilizer of $u$ in $\GG(\Dc)$.  Recall that $\GG(\Dc)$ acts on the
right of words, scanning from right to left.

\begin{lemma}\label{stabilizertransitivity}
Let $p$ be a non-empty pattern and let $p_0$ be the pattern obtained
from $p$ by removing its first letter.  Then $\GG(\Dc)$ acts
transitively on the set of freely irreducible words following $p$ if
and only if it acts transitively on the set of freely irreducible
words following $p_0$ and there is a freely reducible word $w$
following $p$ so that $\St {w_0}$, where $w_0$ is obtained by
removing the first letter of $w$, acts transitively on the set of
freely irreducible words following $p$ with suffix $w_0$.
\end{lemma}
\begin{proof}
Suppose first that $\GG(\Dc)$ acts transitively on the set of freely
irreducible words following $p$.  Then, by ignoring first letters,
it is immediate that it acts transitively on the set of freely
irreducible words following $p_0$. Moreover, if $w$ is any freely
irreducible word following $p$ and $w_0$ is the suffix of $w$
obtained by removing the first letter, then $\GG(\Dc)$ acts
transitively on the set of freely irreducible words following $p$
with suffix $w_0$. But if $xw_0g = yw_0$ then $g\in \St {w_0}$, so
$\St {w_0}$ acts transitively on such words.

For the converse, suppose that $w$ is a word following $p$ as in the
hypothesis and let $w_0$ be obtained by removing the first letter of
$w$.  Let $v$ be a word following the pattern $p$.  Then by
assumption, there exists $g\in \GG(\Dc)$ such that $vg = xw_0$ for
some $x$.  Moreover, $xw_0$ is freely irreducible and follows $p$ by
Lemma~\ref{actpreservespatt}.  Then our assumption on $\St {w_0}$
gives us an element $g'\in St_{w_0}$ with $xw_0g' = w$. Thus $vgg' =
w$. This establishes the transitivity on the pattern $p$.
\end{proof}

Now we turn to the ``Critical Lemma", whose consequences shall be
used over and over again.  It allows us to use the powerful results
of~\cite{free1} concerning the dual to Aleshin's automaton.

\begin{lemma}[Critical Lemma]\label{criticallemma}
 Let $H\leq
\GG(\Dc)$ be the subgroup generated by $\Gc_0,\Gc_1$, $(ab)^\pm$ and
$(bc)^\pm$ .  Let $a_1,\ldots,a_m,b_1,\cdots,b_m\in \{a,b,c\}$,
$x_1,x_m\in ((Q_n\setminus\{a,b,c\})^\pm)^*$, $x_2,\ldots,x_{m-1}\in
(Q_n\setminus\{a,b,c\})^\pm$ and $e_1,\ldots,e_m\in \{\pm 1\}$. Then
there exists $g\in H$ such that
\begin{equation}\label{criticallemmaeq}
x_ma_m^{e_m}x_{m-1}\cdots x_2a_1^{e_1}x_1g=x_mb_m^{e_m}x_{m-1}\cdots
x_2b_1^{e_1}x_1
\end{equation}
\end{lemma}
\begin{proof}
Let $\Gamma$ be the group of transformations of
$(\{a,b,c\}^{\pm})^*$ generated by $\mathsf E_0$, $(ab)^\pm$ and
$(bc)^\pm$ (see Figure~\ref{aleshinE}).  This is the same as the
group generated by the dual automaton to Aleshin's automata by
Proposition~\ref{generatorsforaleshindual}.  By
Theorem~\ref{mainaleshinresult}, $\Ga$ acts transitively on
$\{a,b,c\}^*$.  We basically simulate the action of $\Ga$ using $H$.

More precisely, we show that if $\gamma$ is a generator of $\Ga$ and
$a_m\cdots a_1\gamma = b_m\cdots b_1$, then there exists $g\in H$ so
that \eqref{criticallemmaeq} holds.  The result then follows from
the transitivity of $\Ga$ on $\{a,b,c\}^*$.  If $\ga=(ab)^{\pm}$ or
$(bc)^{\pm}$, then we can take $g=\ga$.  This leaves us with
$\ga=\mathsf E_0$.

Suppose first that that $x_1$ has odd length.  Then
\[x_ma_m^{e_m}x_{m-1}\cdots x_2a_1^{e_1}x_1\Gc_1
=x_ma_m^{e_m}x_{m-1}\cdots x_2a_1^{e_1}\Gc_0x_1\] next we observe
that $x_{i+1}a_i$ switches states in $\Gc$ if and only if $a_i\in
\{a,b\}$ while $a_i$ switches states in $\mathsf E$ if and only if
$a_i\in \{a,b\}$.  Moreover, $\mathsf E_0$, $\mathsf E_1$ act on
letters from $\{a,b,c\}$ as the identity, respectively, as $(ab)$.
On the other hand $\Gc_0$, $\Gc_1$ act on letters from
$\{a,b,c\}^\pm$ as the identity,  respectively, $\overline{(ab)}$.
Thus $x_{i+1}a_i^{e_i}\Gc_j = x_{i+1}(a_i\mathsf E_j)^{e_i}$.
Combining these two observations we see that taking $g=\Gc_1$ gives
the equality \eqref{criticallemmaeq}. If $x_1$ has even length, then
\[x_ma_m^{e_m}x_{m-1}\cdots x_2a_1^{e_1}x_1\Gc_0
=x_ma_m^{e_m}x_{m-1}\cdots x_2a_1^{e_1}\Gc_0x_1\] and the same
argument applies to show that $g=\Gc_0$ does the job.
\end{proof}

Let us state a corollary that we shall use frequently in the sequel.
\begin{corollary}\label{criticalford}
Let $w=d_1^{e_m}a_md_1^{e_{m-1}}\cdots d_1^{e_1}a_1d_1^k$ with
$a_j\in \{a,b\}^\pm$, $e_j\in \{\pm 1$\}, all $j$, and $k\in
\mathbb{Z}$ (we admit the possibility $w=d_1^k$, i.e.\ $m=0$). Then
$a^ew$, $b^ew$ and $c^ew$, where $e\in \{\pm 1\}$ is fixed, are in
the same orbit of $\Gc(\Dc)$. Moreover, if $n>1$, then, for $1<i\leq
n$, $d_i^ew$ is also in this orbit.
\end{corollary}
\begin{proof}
The Critical Lemma immediately applies to show that $a^ew$, $b^ew$
and $c^ew$ are in the same orbit.  Suppose now $n>1$.   First
observe, since $(\{a,b\}^\pm)^*$ is invariant under the action of
$\Fc_0,\Fc_1$:
\[d_i^ew\Fc_0^{n-i+1}=c^ed_{n-i+2}^{e_m}b_md_{n-i+2}^{e_{m-1}}\cdots
d_{n-i+2}^{e_2}b_2d_{n-i+2}^k = c^eu\] with the $b_j\in
\{a,b\}^\pm$. Now by the Critical Lemma, there exist $g_1,g_2\in H$
so that $c^eug_1 = a^eu$ and $c^eug_2 = b^eu$.  Again using that
$(\{a,b\}^\pm)^*$ is invariant under the action of $\Fc_0,\Fc_1$, we
see that one of $a^eu\Fc_0^{-(n-i+1)}$, $b^eu\Fc_0^{-(n-i+1)}$  is
$a^ew$ and the other is $b^ew$. This completes the proof.
\end{proof}

A further property of the subgroup $H$ that we shall need later is
contained in the following lemma.

\begin{lemma}\label{propsofH}
Let $x\in Q^\pm$ and $w\in (Q_n^\pm)^*$.  Suppose that $h\in H$ and
$wxh = uy$ where $y\in Q_n$.   Then there exists $h'\in H$ such that
$wxxh' =uyy$.
\end{lemma}
\begin{proof}
By an easy induction on the length of $g$, it suffices to verify
this for the generators $(ab)^\pm$, $(bc)^\pm$, $\Gc_0$, $\Gc_1$ of
$H$.  Actually, Lemma~\ref{StudyG} shows that $\Gc_1=\Gc_0(ab)^\pm$
so we can omit $\Gc_1$.  If $h=(ab)^\pm,(bc)^\pm$, we can clearly
take $h'=h$.  If $h=\Gc_0$ and $x\in \{a,b\}^\pm$, then we can again
take $h'=h$ since $x$ labels a loop at each vertex of $\Gc$.  If
$x=c^{\pm 1}$, then take $h'=\Gc_1$.
\end{proof}

Another technical lemma that we shall need is the following variant
on Corollary~\ref{criticalford}.
\begin{lemma}\label{hardcorecasebign}
Let $e\in \{\pm 1\}$  and $u\in (\{a,d_1\}^{\pm})^*$. Then $b^eu$
and $c^eu$ are in the same orbit under $\GG(\Dc)$.  Moreover, if
$n>1$ and $1<i\leq n$, then $d_i^eu$ is in the same $\GG(\Dc)$-orbit
as $b^eu$ and $c^eu$.
\end{lemma}
\begin{proof}
Clearly $b^eu(bc)^\pm = c^eu$ and so they are in the same orbit.
Suppose $n>1$ and $1<i\leq n$.     First observe that
$d_i^eu(ac)^{\pm} = d_i^eu'$ where $u'\in(\{c,d_1\}^{\pm})^*$. Then
$d_i^eu'\Fc_0^{n-i+1} = c^eu''$ where
$u''\in(\{d_{n-i+1},d_{n-i+2}\}^{\pm})^*$.  Now $c^eu''(ac)^{\pm} =
a^eu''$ and $c^eu''(bc)^{\pm} = b^eu''$.  Since $\Fc_0\inv$
preserves $\{a^{e},b^{e}\}$, one of the pair
$a^eu''\Fc_0^{-(n-i+1)}$, $b^eu''\Fc_0^{-(n-i+1)}$ is $a^eu'$ and
the other is $b^eu'$. Then $b^eu'(ac)^{\pm} = b^eu$, establishing
that $b^eu$ is in the orbit of $d_i^eu$.
\end{proof}

We now turn to the main technical proposition of the paper, where we
prove transitivity on the patterns.

\begin{proposition}\label{actstransitively}
Let $p\in \{\ast,\ast\inv\}^*$ be a pattern.  Then $\GG(\Dc)$ acts
transitively on the set of freely irreducible words in $(Q_n^\pm)^*$
following the pattern $p$.
\end{proposition}
\begin{proof}
The proof goes by induction on the length of the pattern $p$.
Clearly, the statement is true for the empty pattern.  For patterns
of length $1$, the result holds since the permutation $\delta_0$
acts transitively on $Q_n\inv$ while $\delta_1$ acts transitively on
$Q_n$. Suppose that the proposition is true for all patterns of
length $m$ and suppose that $p$ has length $m+1$.  We shall use
without comment throughout that $\sigma^\pm\in \GG(\Dc)$ whenever
$\sigma\in S_{a,b,c}$.

First we handle the case $p=\ast\ast\cdots \ast$ or
$\ast\inv\ast\inv \cdots \ast\inv$. Let $e\in \{\pm 1\}$.  Consider
$w=d_1^{e(m+1)}$.  By Lemma~\ref{stabilizertransitivity} it suffices
to show that $\St {d_1^{em}}$ acts transitively on
$\{x^ed_1^{em}\mid x\in Q_n\}$.  By Corollary~\ref{criticalford}, we
have that, for $1<i\leq n$, $d_i^ed_1^{em}$ is in the same orbit as
$a^ed_1^{em}$, $b^ed_1^{em}$ and $c^ed_1^{em}$.  Thus, all we need
to do is show that we can change the first letter of $w$ to some
other letter, leaving the rest of $w$ alone.    Now
\[d_1^{e(m+1)}\Fc_0\inv(ac)^\pm = c^{e(m+1)}(ac)^\pm = a^{e(m+1)}\]
Set $i=0$ if $e=1$ and $i=1$ if $e=-1$.  Then: \[a^{e(m+1)}\Fc_i =
\begin{cases}(b^ea^e)^{\frac{m+1}{2}} & m+1\ \text{even}\\
a^e(b^ea^e)^{\frac{m}{2}} & m+1\ \text{odd}\end{cases}\]  We break
the proof into two cases, depending on whether $m+1$ is even or odd.
Suppose that $m+1$ is even.  Then $(b^ea^e)^{\frac{m+1}{2}}(ac)^\pm
= (b^ec^e)^{\frac{m+1}{2}}$ and \[(b^ec^e)^{\frac{m+1}{2}}\Fc_i
=\begin{cases}
(a^ed_1^e)^{\frac{m+1}{2}}& c\ \text{active} \\
(a^ed_1^eb^ed_1^e)^{\frac{m+1}{4}} & c\ \text{inactive and}\ 4\mid m+1\\
(b^ed^e)(a^ed_1^eb^ed_1^e)^{\frac{m-1}{4}} & c\ \text{inactive and}\
4\nmid m+1\end{cases}\] By the Critical Lemma, in all cases we can
find an element $g\in H$ the changing the first letter of
$(b^ec^e)^{\frac{m+1}{2}}\Fc_i$ to, say $c^e$, and leaving the
remaining letters alone. Then, by undoing the previous
transformations (not including $g$), we have managed to change the
first letter, and only the first letter, of $w$ as required.

Suppose now that $m+1$ is odd. Then
$a^e(b^ea^e)^{\frac{m}{2}}(bc)^\pm =a^e(c^ea^e)^{\frac{m}{2}}$. Then
\[a^e(c^ea^e)^{\frac{m}{2}}\Fc_i =\begin{cases}
a^e(d_1^ea^e)^{\frac{m}{2}} & c\ \text{active}\\
a^e(d_1^eb^ed_1^ea^e)^{\frac{m}{4}} & c\ \text{inactive and}\ m/2\
\text{even}\\ b^ed_1^ea^e(d_1^eb^ed_1^ea^e)^{\frac{m-2}{4}}
&\text{inactive and}\ m/2\ \text{odd}\end{cases}\] Again, by the
Critical Lemma, we can change the first letter of
$a^e(c^ea^e)^{\frac{m}{2}}\Fc_i$ to $c^e$ and then undo the previous
transformations with the result that the first letter only of $w$ is
changed.  This finishes the case at hand.

The next case we handle is when $p$ begins with $\ast^e\ast^{-e}$,
where $e\in \{\pm 1\}$.  Let $w=a^ed_1^{-e}\cdots$ be the word
following the pattern $p$ starting with $a^e$ and alternating
$a^{\pm 1}$ with $d_1^{\pm 1}$. Let $w_0=d_1^{-e}\cdots$ be the word
obtained from $w$ by removing the first letter $a^e$.  By
Lemma~\ref{stabilizertransitivity}, it suffices to show that $\St
{w_0}$ acts transitively on the words of the form $x^ew_0$ where
$x\in Q_n\setminus \{d_1\}$.  But this is immediate from
Corollary~\ref{criticalford}.

If the last two letters of the pattern $p$ are $\ast^e\ast^{-e}$,
then the previous case and Lemma~\ref{canreverse} show that all
freely irreducible words following $p$ are in one orbit.  Thus we
are left with the case of a pattern of the form
$\ast^{e_2}\ast^{e_2}\cdots \ast^{e_1}\ast^{e_1}$ where $e_1,e_2\in
\{\pm 1\}$.  We begin with the case where the length $m+1$ of our
pattern $p$ is odd.  Consider $w = a^{e_2}d_1^{e_2}\cdots
d_1^{e_1}a^{e_1}$, the alternating word of length $m+1$ in $a^{\pm
1}$ and $d_1^{\pm 1}$ following the pattern $p$ and starting with
$a^{e_2}$. Let $w_0 =d_1^{e_2}\cdots d_1^{e_1}a^{e_1}$ be the word
obtained by deleting the first $a^{e_2}$ from $w$.  By
Corollary~\ref{criticalford}, $w$ is in the same orbit as
$b^{e_2}w_0$, $c^{e_2}w_0$ and $d_i^{e_2}w_0$, for $1<i\leq n$.  We
thus just need to show that we can find an element of $\GG(\Dc)$
that changes the first letter of $d_1^{e_2}w_0$ and leaves the
suffix $w_0$ alone (or failing this, try and prove transitivity by
some dirty trick). Let $u$ be the word obtained from $w_0$ by
changing the last letter from $a^{e_1}$ to $d_1^{e_1}$; so $u=
d_1^{e_2}\cdots d_1^{e_1}d_1^{e_1}$ where the part in the middle
alternates the letters $a^{\pm 1}$ and $d_1^{\pm 1}$. By induction,
there exists $g\in \GG(\Dc)$ such that $w_0g =u$. Hence,
$d_1^{e_2}w_0g = x^{e_2}u$ where $x\in Q_n$.  Suppose first that
$x\in Q_n\setminus \{d_1\}$.  Choose $y\in Q_n\setminus\{d_1,x\}$.
Then, by Corollary~\ref{criticalford}, we can find $g'\in \GG(\Dc)$
so that $x^{e_2}ug' = y^{e_2}u$. Then $d_1^{e_2}w_0gg'g\inv =
z^{e_2}w_0$ with $z\in Q_n\setminus \{d_1\}$ and we are done. So
suppose instead that $x=d_1$. Then $w_1=x^{e_2}u =
d_1^{e_2}d_1^{e_2}\cdots d_1^{e_1}d^{e_1}$, where the middle part
alternates $a^{\pm 1}$ and $d_1^{\pm 1}$. Thus $w_1$ differs from
$d_1^{e_2}w_0$ only in the last letter. Consider the reversals:
\[(d_1^{e_2}w_0)^{\rho} = a^{e_1}d_1^{e_1}\cdots d_1^{e_2}d_1^{e_2}
= a^{e_1}v\ \text{and}\ w_1^{\rho} = d^{e_1}d_1^{e_1}\cdots
d_1^{e_2}d_1^{e_2} = d^{e_1}v\] Since $d_1^{e_2}w_0$ and $w_1$ are
in the same $\GG(\Dc)$-orbit, Lemma~\ref{canreverse} implies that
$(d_1^{e_2}w_0)^{\rho}=a^{e_1}v$ and $w_1^{\rho}=d_1^{e_1}v$ are in
the same orbit. By Corollary~\ref{criticalford}, $a^{e_1}v$ is in
the same orbit as $b^{e_1}v$, $c^{e_1}v$ and $d_i^{e_1}v$, for
$1<i\leq n$.  Thus $\St {v}$ acts transitively on the set of freely
irreducible words of length $m+1$ following the reverse pattern
$p^{\rho}$ with suffix $v$.  Hence, by
Lemma~\ref{stabilizertransitivity}, $\GG(\Dc)$ acts transitively on
the set of freely irreducible words following the pattern
$p^{\rho}$.  Lemma~\ref{canreverse} then shows that $\GG(\Dc)$ acts
transitively on the set of freely irreducible words following the
pattern $p$.

Our final (and most difficult) case arises when the length $m+1$ of
$p$ is even. Then consider $w=a^{e_2}d_1^{e_2}\cdots
a^{e_1}d_1^{e_1}$ the alternating word of length $m+1$ in $a^{\pm
1}$ and $d_1^{\pm 1}$ following $p$. Let $w_0 =d_1^{e_2}\cdots
a^{e_1}d_1^{e_1}$ be the word obtained by removing the first letter
$a^{e_2}$.  By Corollary~\ref{criticalford}, $w$ is in the same
orbit as a $b^{e_2}w_0$, $c^{e_2}w_0$ and $d_i^{e_2}w_0$, for
$1<i\leq n$. So again we just need to show that we can find an
element of $\GG(\Dc)$ that changes the first letter of
$w'=d_1^{e_2}w_0$ and leaves the suffix $w_0$ alone (or find some
more dirty tricks). Let $u$ be the word obtained from $w_0$ by
changing the last letter from $d_1^{e_1}$ to $a^{e_1}$; so $u=
d_1^{e_2}\cdots a^{e_1}a^{e_1}$ where the part in the middle
alternates $a^{\pm 1}$ and $d_1^{\pm 1}$. By induction, there exists
$g\in \GG(\Dc)$ such that $w_0g =u$.  So $w_1=w'g = d_1^{e_2}w_0g =
x^{e_2}u$ with $x\in Q_n$.  Suppose first that $x\in \{a,b,c\}$. Let
$u=u'a^{e_1}a^{e_1}$. By the Critical Lemma, there exists $h\in H$
such that
\[x^{e_2}u'a^{e_1}h = y^{e_2}u'a^{e_1}\]
where $y\in \{a,b,c\}\setminus \{x\}$.  Then, by
Lemma~\ref{propsofH}, we can find $h'\in H$ such that
\[x^{e_2}uh'=x^{e_2}u'a^{e_1}a^{e_1}h' =
y^{e_2}u'a^{e_1}a^{e_1}\]  Then $w'gh'g\inv$ changes the first
letter of $w'$ and leaves the suffix $w_0$ alone, as required. Thus
we are left with the case that $x\in \{d_1,\ldots, d_n\}$.  Assume
first that $x=d_i$ with $1<i\leq n$.   Then
Lemma~\ref{hardcorecasebign} shows that there is an element $g_0\in
\GG(\Dc)$ with $d_i^{e_2}ug_0=b^{e_2}u$. Then $w'gg_0g\inv$ is a
word ending in $w_0$ but with a different first letter than $w'$ and
we are done.

We turn to the case $x=d_1$, so $w_1 = w'g = d_1^{e_2}u$.  Notice
that $w_1$ and $w'$ differ only in the last letter.  Therefore,
\begin{align*}
w'^{\rho} &=d_1^{e_1}a^{e_1}\cdots d_1^{e_2}d_1^{e_2} = d_1^{e_1}r\\
w_1^{\rho} &=a^{e_1}a^{e_1}\cdots d_1^{e_2}d_1^{e_2} = a^{e_1}r
\end{align*}
where $r$ is an alternating word in $a^{\pm 1}$ and $d_1^{\pm 1}$
starting with $a^{e_1}$. By Lemma~\ref{canreverse}, we have that
$w'^{\rho},w_1^{\rho}$ are in the same $\GG(\Dc)$-orbit, since $w'$
and $w_1$ are in the same orbit. Lemma~\ref{hardcorecasebign} shows
that $b^{e_1}r$, $c^{e_1}r$ and $d_i^er$, for $1<i\leq n$, are all
in the same $\GG(\Dc)$-orbit.  The conclusion that we may draw is
that $\St {r}$ has no singleton orbits on the set $Q_n^{e_1}r$.

   Our goal now is to show
transitivity on freely irreducible words following the reverse
pattern $p^{\rho}$. To do this, consider the alternating word
$t=a^{e_1}d_1^{e_1}\cdots a^{e_2}d_1^{e_2}$ following the pattern
$p^{\rho}$.  Let $t_0 = d_1^{e_1}\cdots a^{e_2}d_1^{e_2}$ be the
word obtained by removing the first letter of $t$. We show $\St
{t_0}$ acts transitively on the set of $Q_n^{e_1}t_0$; this will
give the desired transitivity by Lemma~\ref{stabilizertransitivity}.
By Corollary~\ref{criticalford}, $\{y^{e_1}t_0\mid y\in Q_n\setminus
\{d_1\}\}$ is in a single orbit. So it suffices to show that $\St
{t_0}$ has no singleton orbits on $Q_n^{e_1}t_0$.  By induction,
there exists $g'\in \GG(\Dc)$ such that $t_0g'=r$.  Hence the action
of $\St {t_0}$ on $Q_n^{e_1}t_0$ is conjugate to the action of $\St
{r}$ on $Q_n^{e_1}r$ --- but the latter was already shown to have no
singleton orbits.  This establishes the transitivity for $p^{\rho}$
and hence, by Lemma~\ref{canreverse}, for $p$.  This completes the
proof of the proposition.
\end{proof}

As a consequence of Propositions~\ref{actstransitively} and
\ref{actsnontrivlevel1} and in light of Corollary \ref{freeness}, we
have proven our main theorem:

\begin{theorem}\label{ourfreetheorem}
If $\Ac$ belongs to the family $\mathfrak F$, then the states of
$\Ac$ freely generate a free group.  In particular for every even
number $n\geq 4$, there is an $n$-state connected automaton over a
binary alphabet, whose states freely generate a free group of rank
$n$.
\end{theorem}

\subsection{Free products of cyclic groups}
We obtain the result for the family $\mathfrak F'$ from the result
for the family $\mathfrak F$ using a straightforward fact from
combinatorial group theory, a proof of which can be found
in~\cite[Lemma 6.5]{free2}.

\begin{lemma}\label{freetocyclic}  Suppose that $G$ is a group
generated by elements $g_0,g_1,\ldots, g_k$, $k\geq 1$, satisfying
$g_i^2=1$, $0\leq i\leq k$.  Let $H$ be the subgroup generated by
$h_1,\ldots, h_k$, where $h_i = g_0g_i$, $1\leq i\leq k$.  Then $G$
is a free product of $k+1$ cyclic groups of order two, freely
generated by $g_0,\ldots,g_k$, if and only if $H$ is a free group of
rank $k$, freely generated by $h_1,\ldots, h_k$.
\end{lemma}

Theorem~\ref{ourfreetheorem}, in light of
Lemmas~\ref{secondfamilydetails} and \ref{freetocyclic}, establishes
the following result:

\begin{theorem}
Let $\Ac=(Q_n,\{0,1\},\lambda,\delta) \in \mathfrak F$, with $n\geq
1$, and let $(01)[\Ac]\in \mathfrak F'$ be the corresponding
automaton. Then the group generated by the free monoid automorphism
$(01)^\ast$ and the elements $(01)^*\Ac_q=(01)[\Ac]_q$, with $q\in
Q_n$, is a free product of $n+4$ cyclic groups of order two, freely
generated by $(01)^*$ and the $n+3$ states of $(01)[\Ac]_q$.  In
particular, the states of each member of $\mathfrak F'$ freely
generate a free product of cyclic groups of order two.
\end{theorem}

\bibliographystyle{amsplain}

\end{document}